\input amstex
\documentstyle{amsppt}
\magnification=\magstep0
\define\cc{\Bbb C}

\define\r{\Bbb R}

\define\N{\Bbb N}

\define\A{\Cal A}
\define\h{\Cal D}

\define\m{\Cal M}

\define\f{\Cal S}

\define\n{\Cal F}

\define\la{\lambda}
\define\om{\omega}
\define\Om{\Omega}
\define\e{\varepsilon}
\define\va{\varphi }
\define\CB#1{\Cal C_b(#1)}
\define\st{\subset }
\define\al{\alpha}
 \topmatter
 \title
Existence of bounded uniformly continuous mild solutions on $\Bbb{R}$ of evolution equations and some applications.
\endtitle
\endkeywords
 \author
 Bolis Basit and Hans G\"{u}nzler
\endauthor
\abstract{ We prove that  there is $x_{\phi}\in X$ for which (*)$\frac{d u(t)}{dt}= A u(t) +  \phi (t) $, $u(0)=x$ has on $\r$ a mild solution $u\in C_{ub} (\r,X)$ (that is bounded  and uniformly continuous)  with $u(0)=x_{\phi}$, where $A$ is the generator of a holomorphic $C_0$-semigroup $(T(t))_{t\ge 0}$
   on ${X}$ with sup $_{t\ge 0} \,||T(t)|| < \infty$,  $\phi\in L^{\infty} (\r,{X})$ and  $i\,sp (\phi)\cap \sigma (A)=\emptyset$.   As a consequence it is shown that if $\n$ is the  space of almost periodic $AP$, almost automorphic $AA$,
          bounded Levitan almost periodic $LAP_b$, certain classes of  recurrent functions $REC_b$ and $\phi \in  L^{\infty} (\r,{X})$ such that $M_h \phi:=(1/h)\int_0^h \phi (\cdot+s)\, ds \in \n$  for each $h >0$, then $u\in \n\cap C_{ub}$. These results seem new and generalize and strengthen several recent Theorems.}
 \endabstract
\endtopmatter
\rightheadtext{ admissibility} \leftheadtext{ Basit and
 G\"unzler}
 \TagsOnRight

\document
\baselineskip=22pt

\head{\S 1.  Introduction,     Definitions and Notation}\endhead

In this paper
 we study solutions of the inhomogeneous  abstract Cauchy problem

(1.1)\qquad $\frac{d u(t)}{dt}= A u(t) +  \phi (t) $, $u(t)\in X$,  $t\in  {J}\in\{\r_+,\r\}$,

(1.2)\qquad $u(0)=x$,

\noindent  where $A: D(A)\to X$ is the generator of a
 $C_0$-semigroup $(T(t))_{t\ge 0}$
on the complex Banach space ${X}$ and $\phi\in L^1_{loc} ({J},{X})$, $x\in X $.

\noindent By [15, Corollary 2.5, p. 5], it follows that $D(A)$ is dense in $X$ and $A$ is a closed linear operator.

\noindent By a $classical\,\, solution$ of (1.1) we mean a function $u: J\to D(A)$ such that $u \in C^1 (J,X)$ and (1.1) is satisfied.

\noindent  By a $mild\,\, solution$ of (1.1), (1.2) we mean a $\omega  \in
C(J,X)$  with  $\int _ 0 ^ t  \omega(s)ds  \in  D(A)$ for  $t \in J$ and

(1.3)\qquad $\om(t)  =  x  + A \int_0^t  \omega(s)\, ds  + \int_0^t \phi(s)ds$, $t  \in J $.

\noindent For $ J = \r_+$  this is the usual definition [1, p. 120].

A classical solution is always a mild solution (see [1, (3.1), p. 110], for
$J = \r_+$, $\phi = 0$).
Conversely, a mild solution with $\omega  \in  C^1(J,X)$ and $\phi \in C(J,X)$ is
a classical solution (as in [1, Proposition 3.1.15]).

With translation and the case $J = \r_+ $ [1, Proposition 3.1.16] one can show,
for $J, A, T$, $\phi$ as in (1.1) and after (1.2)

\proclaim {Lemma 1.1}  $\omega : J \to X $ is  a mild solution of (1.1) if and only if

(1.4) \qquad $\om(t)= T(t-t_0)\om(t_0)+\int_{t_0}^t T(t-s)\phi(s)\,ds$,\,\,\,$t\ge t_0$,\,\,\, $t_0\in J$.
\endproclaim

 When $J=\r$ mild solutions of (1.1)  may not exist. In this paper we  establish  the existence of a bounded  uniformly continuous  mild solution  on $\r$ to (1.1), (1.2)   for  $\phi\in L^{\infty} (\r,X)$,  $x =x(\phi)\in X$ and $i\,sp(\phi)\cap \sigma (A) =\emptyset$ when $A$ is the generator of  a holomorphic $C_0$-semigroup $(T(t))$ satisfying sup $_{t \in \r_+}||T(t)|| <\infty$.
(Theorem 3.5). In Theorem 3.6, we show that all mild solutions of (1.1), (1.2) are  bounded  and uniformly continuous when $J=\r_+$ and $\phi$ and $A$ as in Theorem 3.5. Theorem 3.6(i) should be
compared with the "Non-resonance Theorem 5.6.5 of [1]" :
the result of [1] is  more general since the half-line spectrum is used ($ \st$  Beurling spectrum). However, our main result (Theorem 3.5) can not be deduced from Theorem 5.6.5 of [1]. Moreover, in the important cases of almost periodic, almost automorphic and recurrent functions,   the half-line spectrum  coincides with  the  Beurling spectrum (see [8, Example 3.8] and [10, Corollary 5.2]).

By mild admissibility  of $\n$ with respect to equation (1.1)  we mean that for every $\phi\in \n$, where $\n \st L^{\infty}(\r,X)$ is a subspace  with (4.2)- (4.5) below equations (1.1), (1.2) have a unique mild solution  on $\r$ which belongs to $\n$ (see [14, Definition 3.2,  p. 248]). In the case $\n \st C_{ub} (\r,X)$, several methods have been used to prove admissibility for (1.1), (1.2) (see  [16], [4], [14], [12] and references therein).  The method of sums of commuting operators  was introduced  in [14], and
 then extended to the case $\n= AA(\r,X)$ in [12], to prove admissibility for (1.1). As consequences of Theorem 3.5  we generalize and strengthen some recent results on admissibility. Namely, we prove that if $\phi \in \m\n \cap  L^{\infty} (\r,X)$  (see (4.1)) and  $i\,sp(\phi)\cap \sigma (A) =\emptyset$, then there is a mild solution  $u_{\phi}\in \n \cap C_{ub} (\r,X)$ of (1.1), (1.2). For example if $\phi \in \m AP \cap L^{\infty} (\r,X)$, then $u_{\phi}\in AP$ or if $\phi \in \m AA \cap L^{\infty} (\r,X)$, then $u_{\phi}\in AA \cap C_{ub}$.

Further notation and definitions.
$\r_+ = [0,\infty)$, $L(X) =\{B : X\to X,   B $ linear bounded $\}$ with norm
$||B||$,  $\f (\r)$ contains  Schwartz's complex valued rapidly decreasing $C^{\infty }$-functions,
$C_{ub}(\r,X) = \{ f : \r\to X : f$  bounded uniformly  continuous $\}$, $AP = AP(\r,X)$ almost periodic functions, $AA = AA(\r,X)$ almost automorphic functions [19], $BAA = BAA(\r,X)$ Bochner almost automorphic functions  [21, p. 66], [12]; for
$f \in  L^1_{loc}(J,X)$  $(Pf)(t) : =$
Bochner integral  $\int^t_0 f(s)\, ds$, $\widehat {f}(\lambda) = L^1$-Fourier transform
 $\int_{\r} f(t) e^{- i \lambda\, t}\, dt$,  $f_a(t) =$  translate $f(a+t)$ where defined, $a$ real, $sp
=$  Beurling spectrum (2.3), Proposition 2.3.

\head{\S 2. Preliminaries}\endhead

In this section we collect some  lemmas and propositions needed in the sequel.

\proclaim{Proposition 2.1}  Let  $F\in L^1(\r,L(X))$ and  $\phi\in L^{\infty}(\r,X)$ respectively  $F\in L^1(\r,X)$ and  $\phi\in L^{\infty}(\r,\cc)$ respectively  $F\in L^1(\r,\cc)$ and  $\phi\in L^{\infty}(\r,X)$.

(i) $F*\phi (t):= \int_{-\infty}^{\infty}F(s)\phi (t-s)\, ds$

\noindent exists as a Bochner integral for $t\in \r$ and
$F*\phi\in C_{ub}(\r,X)$.

(ii) If moreover $ P\phi\in L^{\infty}(\r,X)$ respectively   $P\phi\in L^{\infty}(\r,\cc)$, then

(2.1) \qquad $F*(P\phi)\in C^1(\r,X)$, $(F*(P\phi))'=F*\phi$,

(2.2)\qquad  $F*(P\phi)= P(F*\phi)+(F*(P\phi))(0)$.

\endproclaim

\demo{Proof} (i) The integrand of (i) is measurable (approximate $F$ and $\phi$ by
continuous functions a.e),  it is dominated  by  $||\phi||_{\infty}  |F|   \in L^1(\r,\r)$, $|F|(t) : =
||F(t)||$, so this integrand $ \in L^1(\r,L(X))$ for any fixed  $t \in\r$.
To $F$ there is a sequence $(H_n)$ of $L(X)$-valued step-functions with
$||F - H_n|| _{L^1}  \to 0$. It follows $||F*\phi  -  H_n*\phi||_{\infty}  \to 0 $, so
with  the  $H_n*\phi$ also $F*\phi$ is bounded and uniformly continuous. Similarly in the other cases.

(ii) By part (i) we conclude that $F*P\phi,  F*\phi \in C_{ub}(\r,X)$; [1, Proposition 1.2.2] and the Lebesgue convergence theorem give existence
      and
 $(F*P\phi)'
=  F*\phi, \in C_{ub}(\r,X)$, so  $(F*P\phi)  \in C^1(\r,X) $. It follows $F*P\phi=P (F*\phi) + F*P\phi (0)$. \P
\enddemo

In the following $sp$ denotes the Beurling spectrum, $sp (\phi) = sp_{\{0|R\}} (\phi)$ as defined for example in
 [8, (3.2), (3.3)] case $S=\phi\in L^{\infty} (\r,X)$, $V= L^1 (\r,\cc)$, $\A= \{0|\r\}$:

(2.3) \qquad $sp_{\{0|\r\}} (\phi):= \{\om\in \r: f\in L^1 (\r,\cc), \phi*f =0 $ imply $\widehat {f}(\om)=0\}$.

$sp_B$ is defined in [1, p. 321], $sp_C$ is the Carleman spectrum
[1, p.293/317].

\proclaim {Proposition 2.2}  If  $\phi\in L^{\infty} (\r,X)$, then

(2.4) \qquad  $sp (\phi)  =  sp_{\{0|R\}} (\phi)  =  sp_B (\phi)  =  sp_C (\phi) $.

 See also [8, (3.3), (3.14)].
\endproclaim
\demo{Proof} $sp (\phi)  \st  sp_B \phi $:   Assuming $\la \in sp (\phi)$ and  $h\in L^1 (\r,\cc)$ with $\widehat {h} (\la)\not = 0$,  we conclude $\phi*h \not =0$. This  implies $\la \in sp_{B} (\phi)$.

 $sp_{B} (\phi)\st sp (\phi)$: Assume $\la \in sp_{B} (\phi)$ but $\la \not\in sp(\phi)$. Then there is $h\in L^1 (\r,\cc)$ with $\widehat {h} (\la)\not = 0$ and $\phi*h =0$.  With Wiener's
inversion theorem [11, Proposition 1.1.5 (b), p. 22], there is $h_{\la}\in L^1$ such that $k=h*h_{\la}$ satisfies
 $\widehat{k}=1$ on some neighbourhood $V = (\la-\e,\la+\e)$ of $\la$. Now let $g\in L^1 (\r,\cc)$ be such that supp$\, \widehat {g}\st V$. Then $k*g=g$ and $0= \phi*h = \phi*k = \phi *(k*g)= \phi*g$. This implies $\la\not \in sp_{B}(\phi)$ by [1, p. 321]. This is a contradiction.

$ sp_B (\phi) = sp_C(\phi)$:  See  [1, Proposition  4.8.4, p. 321]. \P
\enddemo

\proclaim{Lemma 2.3}  Let  $\phi \in L^{\infty}(\r,X)$.

(i) If  $P\phi \in L^{\infty}(\r,X)$, then
 $ sp(\phi)\st sp(P\phi)\st sp(\phi)\cup \{0\}$.

(ii) If  $F  \in L^1(\r,\cc)$ or $F\in L^1 (\r,L(X))$, then $ sp(F*\phi)  \st   sp (\phi) \cap $ supp $\widehat{F}$.

(iii) If  $f  \in L^1(\r,\cc)$, then $ sp(\phi - \phi*f)  \st   sp (\phi) \setminus U$, where $U$ is the interior of  the set $\{ \la\in \r : \widehat{f}(\la) = 1 \}$.

 See also [3, Theorem 4.1.4], [14, Theorem 2.1], [17, Proposition 0.4].
\endproclaim

\demo{Proof} (i) If  $f  \in L^1(\r,\cc)$ with  $(P\phi)*f = 0 $,  then  (2.4)  yields
         $\phi*f = 0$ and  so the first inclusion  follows from (2.3((see also [17, Proposition 04 (iv), p. 20]). The second inclusion follows by [9,  Proposition 1.1. (f)] valid also for $
X$-valued functions.

(ii)  $F \in L^1(\r,L(X))$: to $\lambda \not \in sp  (\phi)$ exists $f \in L^1(\r,\cc)$ with $\phi*f=0$,
     $\widehat{f}(\lambda)=1$, then  $(F*\phi)*f = F*(\phi*f) = 0$, so $\lambda \not \in
     sp(F*\phi)$, yielding  $ sp(F*\phi)  \st  sp (\phi) $. (The associativity of convolution used here can be shown as in [1, Proposition 1.3.1, p. 22]).
     If  $\lambda \not \in supp \widehat{F}$, there exists $f \in \f(\r)$ with $(\text {supp\,} \widehat{f}) \cap\text {supp\,} \widehat{F}
     =  \emptyset$ and $\widehat{f}(\lambda)=1$, then $\widehat {F*f} = \widehat {F}\widehat {f} = 0$, and so $F*f = 0$;  this gives
      $ 0 =  (F*f)*\phi = F*(f*\phi) = F*(\phi*f) = (F*\phi)*f$,  and so $\lambda
      \not \in  sp(F*\phi)$, yielding   $sp(F*\phi)  \st \text{ supp\,}\widehat {F}$.
      (Existence of  $F*f,  F*f  \in  L^1(\r,L(X))$, and  $\widehat {F*f} = \widehat {F}\widehat {f}$ follows
       with Fubini-Tonelli similar as for Proposition 2.1(i).)

The proof for the case $F  \in L^1(\r,\cc)$ is similar.

(iii) If $\la\in sp (\phi) \cap U$, then there is $f_{\la}\in \f(\r)$ with supp$\, \widehat{f_{\la}} \st U$ such that $\widehat{f_{\la}} (\la)\not = 0$ and $f*f_{\la}=f_{\la}$. It follows $\la \not \in sp(\phi-\phi*f)$.
 \P
\enddemo

In the following if $U, V \st \r$, then $d(U,V):=$ inf $_{u\in U, v\in V} |u-v|$. If $U, V$ are compact and $U\cap V =\emptyset$, then $d(U,V) > 0$.

\proclaim{Lemma 2.4}  Let  $\psi\in C_b(\r,X)$.  Then there exists a sequence of $X$-valued trigonometric polynomials $\pi_n $ such that

(i) \qquad sup$_{n\in \N} ||\pi_n||_{\infty} < \infty$.

(ii) $\pi_n (t)\to \psi (t)$ as $n\to \infty$ locally uniformly on $\r$.

(iii) If the Beurling spectrum $sp(\psi)\cap [\al,\beta] =\emptyset$, then there exists $\delta > 0$  such that the Fourier exponents of $\pi_n$ can be selected from $\r\setminus (\al-\delta,\beta+\delta)$.

(iv) If $K_1, K_2 $ are compact subsets of $(\al,\beta)$, where $ K_1 =sp(\psi) $ is  the Beurling spectrum of $\psi$ and  $K_1 \cap K_2 =\emptyset$ and $\eta = d(K_1,K_2)$,   then there exist  $W$,  $\va$ so that the Fourier exponents of $\pi_n$
                    can be selected from $ W$, where $W$ is an open set containing $K_1$, and  $\va\in \h(\r)$ with $\va = 1$ on W  and
               $d(K_2,W) \ge  d(K_2,\text{\,supp\,} \va) \ge  \eta/2$.
\endproclaim

\demo{Proof} (i), (ii): Similar to the proof of Lemma 3.1  (3.1) of [4] with [1, Theorem 4.2.19].

(iii) Choose $(\al_1,\beta_1)\st \r$ such that $[\al,\beta]\st (\al_1,\beta_1)$ and $sp(\psi) \st \r\setminus (\al_1,\beta_1)$. Let $f\in \f(\r)$ be such that  supp $\widehat{f} \st (\al_1,\beta_1)$ and $\widehat{f}=1$ on  $[\al-\delta,\beta+\delta]$ for some $\delta >0$. Let $(T_n)$ satisfy (i), (ii).
Then  $\pi_n= T_n- T_n*f$ satisfy (i), (ii) since by Lemma 2.3 (ii) $sp(\psi*f) \st (\r\setminus (\al_1,\beta_1))\cap (\al_1,\beta_1)=\emptyset$ implies $\psi*f=0$  by [1, Theorem 4.8.2 (a)] and Proposition 2.2, and the Fourier exponents of $\pi_n$ belong to $\r\setminus (\al-\delta,\beta +\delta)$.

(iv)
 Let $\Gamma =(I_{t})_{t\in K_1}$, $\Om= \cup_{t\in K_1} I_{t}$ be a system of open intervals
with the properties

\qquad $I_{t}= (t-\delta, t+\delta)$,
   $I_{t}\st (\al,\beta)$  and $ 0 < \delta < \eta/2$ for each $ t\in  K_1$.

 \noindent With a partition of unity     [18, Theorem 6.20, p. 147] there is a  sequence $(\psi_{i})\st \h(\Om)$ with $\psi_{i}$ supported in some  $I_{t_i}\in \Gamma$,  ${\psi_{i}} \ge 0$  and  $\Sigma_{i=1}^{\infty}{ \psi_{i}}(t)=1$, $t\in \Om$. Moreover, for $K_1$ there is   a positive integer $m$  and an open set $W\supset K_1$ such that

(2.5)\qquad  $\va :=\Sigma_{i=1}^m {\psi_{{i}}}(t)=1$ for each $t\in W$.

\noindent  Since  supp $\psi_i  \st  I_{t_i}$,  $d(K_2,W)\ge
                           d(K_2,\text{\,supp\,} \va)  \ge \eta/2 $.
Now consider the system $\Gamma_W:= (J_t)_{t\in K_1}$, where  $J_t= I_t \cap W $ for each $t\in K_1$. As above there is an open set $U\supset K_1$ and $\va_W \in \h (\r)$ such that $\va_W =1$ on $U$ and  supp $\va_W \st W$. It follows  $\va_W= \widehat{f_W}$ with $f_W\in \f (\r)$. Let  $(T_n)$ be a sequence of $X$-valued trigonometric polynomials satisfying (i), (ii).
Take $\pi_n = T_n *f_W$, $n\in \N$. Then $sp(\pi_n)\st W$, Lemma 2.3 (iii)  gives $sp (\psi- \psi*f_W)= \emptyset$  and so $\psi= \psi* f_W$ as in (iii).
\P
\enddemo

\proclaim{Proposition 2.5} If $f,\widehat{f} \in L^1 (\r,X)\cap C_b(\r,X)$, then ([1, Theorem 1.8.1 d), p. 45])

$f(t) = (1/2\pi)\int_{\r} \widehat{f} (\la) e^{i\la t}\, d\la$, $t\in \r$.
\endproclaim

\head{\S 3.  Holomorphic $C_0$-semigroups}\endhead

In this section we study (1.1) when $A$ is the generator  of a holomorphic  $C_0$-semigroup $(T(t))$in the sense  [1, Definition 3.7.1]. By [1,  Corollary 3.7.18], it follows that there is $ a >0$  such that

(3.1)  $\{is, \, s\in \r: \, |s| > a \} \st   \rho (A)$ and

(3.2)  sup $_{|s| > a }\,  ||s\, R(is,A)|| < \infty$

\proclaim{Lemma 3.1} Let $b >0$, $m\in \N$. There are $f_k\in C^{m} ([0,b])$, $k=0, 1, \cdots, n$ such that

$f_k^{(k)} (b)=1$, $f_k^{(j)}(b)=0$, $k\not =j$, $k, j= 0,1,\cdots, m$,

 $f_k^{(j)}(0)=0$, $k, j= 0,1,2,\cdots, m $,  ($f^{(0)}=f$).
\endproclaim

\demo{Proof} The case $m=2: $ Let $f(t)= (\pi/2) \sin (\pi (t/b)^3/2)$, $t\in  [0,b]$. Then

$f_0 (t) =\sin f(t)$, $f_1 (t) =(t-b) f_0(t)$, $f_2 (t) =(1/2)(t-b)^2 f_0(t)$

\noindent satisfy all the requirements of the statement.

The general case follows by   the Lagrange Interpolation Theorem  ([2,  p. 395]). \P
\enddemo

\proclaim{Lemma 3.2}  Let $A: D(A)\to X$ be the generator of a holomorphic $C_0$-semigroup $(T(t))$.

(i) There exists $a > 0$ such that
$R(i\la, A)= (i\la I-A)^{-1}$ is infinitely differentiable  on $\r\setminus (-b,b)$ for each $b > a$. Moreover,
$R(i\la, A)|(\r\setminus (-b,b))$ can be extended to a function  $H_m\in C^m (\r,L(X))$ for any $m\in \N$ with $  H_m^{(j)}(0)=0$, $0\le j \le m$, $H_m(t)X \st  D(A)$
                 for $t  \in\r $.

(ii) $H_m^{(k)}\in L^1 (\r,L(X))$, \,\,\, $1 \le k \le m$.
\endproclaim

\demo{Proof} (i) By (3.1) $ \sigma (A)\cap i\,\r \st i\,[-a,a]$ for some $a > 0$.  It follows $R(i\la, A)= (i\la I-A)^{-1}$ is infinitely differentiable on $\r\setminus (-b,b)$ for each $b > a$.
 Fix one such $b$ and let $f_k (t)$ satisfy the conditions of Lemma 2.1 and $g_k (-t) = f_k(t)$, $t\in [0,b]$, $k=0, 1, \cdots, m$. Set

 $F(t)= f_0 (t)R(ib, A)+f_1 (t)R'(ib, A)+\cdots + f_m(t)R^{(m)}(ib, A)$, $t\in [0,b]$,

 $F(t)= g_0 (t)R(-ib, A)+g_1 (t)R'(-ib, A)+\cdots + g_n(t)R^{(m)}(-ib, A)$, $t\in [-b,0]$.

\noindent Then
 $F\in C^{m}([-b,b],L(X))$.
  Set

 (3.3)\qquad $H_m (\la)= R(i\la, A)$ on $\r\setminus (-b,b)$ and $H_m=  F $ on $[-b,b]$.

\noindent   Then $H_m\in C^{m}(\r,L(X))$. Since  $R(\lambda,A)X  \st   D(A)$ by definition ([1, p. 462]) and $D(A)$ is
           linear, $H_m(\r)X \st  D(A)$ with Lemma 3.1.

(ii)  By(3.2) there  is  $M = M(m) >0$ with

(3.4) \qquad $||(R(i\la, A))^k|| \le M/|\la|^k $ on $\r\setminus (-b,b)$,\,\,\, $1 \le k \le m$.

  \noindent Since ( [1, p.
     463, Corollary B.3])

(3.5)\qquad  $H_m^{(k)}(\lambda) =  k!(-i)^k (R(i\lambda,A))^{k+1}$ if $b \le |\lambda|$,\,\,\, $ 1\le k
 \le m$,

 \noindent with (3.4),(3.5) and part (i) one has $ H_m^{(k)}  \in  L^1(\r,L(X))$, $ 1\le k
 \le m$. \P
\enddemo

\proclaim{Proposition 3.3}  Let $m \ge 2$ and $H_m$ be as in (3.3).

 (i) The function $G_m$ defined by

(3.6)\qquad $G_m(t):= (1/2\pi) \int_{\r} H_m(\la) e^{i\la t}\, d\la$,

\noindent exists as an  improper Riemann integral  for each $t\not = 0$.

(ii)\qquad  $t^k G_m \in L^1 (\r,L(X))$, $k=1, 2,\cdots, m-2$.
\endproclaim

\demo{Proof}     With Lemma 3.2 (ii),  one has $H_m', H_m''  \in  L^1(\r,L(X))$.
    $H_m$ being continuous on $\r$, for $t \not =0$ and $T>0$ one can use partial
integration and gets

     $\int_0 ^ T  H_m(\lambda)  e^{i \lambda t}\, d\lambda
     =(1/it) H_m(T)e^{itT}  -(1/(it)) \int_0 ^ T   H_m'(\lambda) e^{i \lambda t}\, d\lambda$,

 \noindent    similarly for the integral from  $-t$ to $0$.
 This shows that  $G_m(t)$  given by (3.6)  exists as an improper Riemann integral
      and is given by

        $(i/(2t\pi)) ( \int_{ - \infty} ^{ \infty}  H_m'(\lambda)
       e^{i \lambda t}\, d\lambda $, $t\not=
       0 $.

 \noindent   Induction gives with Lemma 3.2 (ii), for $1\le k \le m$,

(3.7) \qquad $G_m(t)=(1/2i\pi)\, (i/t)^k \int_{\r} H_m^{(k)} (\la) e^{i\la t}\, d\la$, $t \not= 0$.

\noindent Hence $G_m$ is continuous  and $||G_m(t)||\le 1/(2\pi \,|t|^k) \int_{\r}||H_m ^{(k)}(\la) ||\, d\la$ in each point $t\not = 0$.

(ii) This follows by (3.7) using Lemma 3.2 (ii). \P
\enddemo

\proclaim{Lemma 3.4} (i) $I (t)=\int_b^{\infty} R(i\la,A) \cos \la t\, d\la$ exists
as improper  Riemann integral for each $t> 0$. Moreover, $||I(t)||= O (|\ln t|)$ for all $ 0 < t < \pi/2b$.

(ii) $J (t)=\int_b^{\infty} R(i\la,A) \sin \la t\, d\la$ exists
as improper  Riemann integral for each $t> 0$. Moreover, $J(t)$  is bounded for all $ t >0$.

(iii) For $G_m$ of Proposition 3.3 and $m \ge 2$, one has $||G_m(t)||= O(|\ln |t||)$ in  some neighbourhood of $0$, $G_m\in L^1 (\r,L(X))$ (Bochner-Lebesgue integrable functions).

(iv)  $\widehat{G_m}= H_m$ if $m \ge 3$ for the $
H_m$ of Lemma 3.2.

(v) $G_m*x=0$ for each $x\in X$, $m\ge 2$.
\endproclaim

\demo{Proof} (i) For  $  t \ge \pi/2b$,  $I(t)$ exists similarly as in Proposition 3.3. For  $ 0 < t < \pi/2b$, we have  $I (t)= (1/t) \int_{bt}^{\infty} R(i\la/t, A) \cos \la\, d\la=(1/t) [a_0 (t)+ \sum _{k=1}^{\infty} a_k (t)]$,   where

 $a_0 (t)=\int_{bt}^{\pi/2}R(i\la/t, A) \cos \la  \, d\la$,

  $a_k (t)=\int_{\pi(1/2 +2 (k-1))}^{\pi (1/2 +2 k)} R(i\la/t, A) \cos \la  \, d\la $
 $= \int_{0}^{\pi}(B_k(\la,t)-A_k(\la,t)) \sin \la  \, d\la $, with

$B_k(\la,t)=R(i(\la+\pi(3/2 +2 (k-1)))/t,A)$,

 $A_k(\la,t)=R(i(\la+\pi(1/2 +2 (k-1)))/t,A)$.

Since $R(\mu,A )-R(\la,A)= (\la-\mu)R(\mu,A)R(\la,A)$ for any $\la, \mu \in \rho (A)$  [1, B.4, p. 464],  we have
$B_k(\la,t)-A_k(\la,t)=$

$ - (i\pi/t) R(i(\la+\pi(3/2 +2 (k-1)))/t,A) R(i(\la+\pi(1/2 +2 (k-1)))/t,A)$.

\noindent Using (3.4) we  conclude  that
$||a_k (t)||= O  (\frac{t} {k^2})$ for all $t > 0$, $k \in \N$. It follows
that series $(1/t)\sum_{k=1}^{\infty} a_k$ is convergent and bounded on $(0,\infty)$. Again by (3.4),

$(1/t)||a_0|| \le  M\int_{bt}^{\pi/2} \frac {\cos \la}{\la}  \, d\la \le M\int_{bt}^{\pi/2} \frac {1}{\la}  \, d\la =  M\ln (\pi/2)-\ln(bt) = O(|\ln t|)$.

(ii) we have  $J (t)= (1/t) \int_{bt}^{\infty} R(i\la/t, A) \sin \la\, d\la=(1/t) [b_0 (t)+ \sum _{k=1}^{\infty} b_k (t)]$,   where
 $b_0 (t)=\int_{bt}^{\pi}R(i\la/t, A) \sin \la  \, d\la$,\,\,\,
$b_k (t)=\int_{\pi(1+2 (k-1))}^{\pi(1+2 k)} (R(i\la/t, A) \sin \la \, d\la $.

\noindent Since $\frac{\sin\la}{\la}$ is bounded on $\r\setminus\{0\}$, we conclude that $b_0(t)/t$ is bounded for all $t> 0$.  Similarly as in part (i), we  conclude that the series $(1/t)\sum_{k=1}^{\infty} b_k (t)$ is convergent for all $t >0$ and bounded on $(0,\infty)$.

(iii) Follows  as in parts (i), (ii) and Proposition 3.3 noting that

$2\pi G_m(t)=   \int_{-b}^{b}H_m(\la) e^{i\la t}\, d\la+  \int_{-\infty}^{-b} R(i\la, A) e^{i\la t}\, d\la + \int_{b}^{\infty} R(i\la, A) e^{i\la t}\, d\la)= $

  $[I(b,t)+\int_{b}^{\infty} (R(-i\la, A)+R(i\la, A)) \cos \la t + (R(i\la, A)-R(-i\la, A)) \sin \la t)\,d\la]$,

  \noindent where  $I(b,t)=\int_{-b}^{b} H_m(\la) e^{i\la t}\, d\la$.

(iv) Since $H_m', t G_m\in L^1(\r,L(X))\cap C(\r,L(X))$ by Proposition 3.3(ii) and

$-i (tG_m)(t)= (1/2\pi)\int_{\r} e^{i\la t}  H_m'(\la)\, d\la $,

\noindent one gets existence of $\widehat{G_m}'$ and by Proposition 2.5 for $\la \in \r$

$H_m'(\la)=-i \int_{\r} e^{-i\la t}  tG_n(t)\, dt =  \widehat {G_m}' (\la)$.

\noindent This implies $H_m =\widehat{G_m}$.

(v) $G_m*x (t)= \int_{\r} G_m(t)x\, dt= \widehat{G_m x}(0)= \widehat{G_m}(0)x= H_m(0)x=0$, with Lemma 3.4, Proposition 2.1 and Lemma 3.2.
 \P
\enddemo

In the following we use:

For any  $G\in L^1 (\r,L(X))$, $\la\in \r$ with $\widehat{G}(\la)= R(i\la,A)$  and $x \in X$ one has  $Gx \in L^1(\r,X)$  and
   $\widehat{(Gx)}  = \widehat{(G)}x$ [1, Proposition 1.1.6], so for
     $e_{\lambda}(t) : = e^{i \lambda t} x$, the following LB-integrals exist and

   $(G*e_{\lambda})(t) = e^{i\lambda t} \int_{\r} G(s)x  e^{- i \lambda s} ds =
   e^{i \lambda t}\widehat {(Gx)}(\lambda) = e^{i \lambda t}\widehat{G}(\lambda)x =
  e^{i \lambda t} R(i \lambda,A) x$,

  \qquad \qquad  \qquad \qquad  \qquad \qquad  $t\in \r$, with $R(i \lambda.A)x \in D(A)$.

\noindent  So for  $x\in X$ and the $\lambda $ and the $G$ above, the equation

(3.8) \qquad
$u'(t)= Au (t)+ x e^{i\la t}$,  $t\in \r$, $u(0)= R(i\la,A)x$

\noindent has a  classical solution $u(t) = e^{i \lambda t} R(i \lambda,A)x =
(G*e_{\lambda})(t)$, $t\in\r $.

Let $A$, $a, b$ be as in Lemma 3.2, $\delta =b-a$ and  $c = b+3\delta$.  Choose $f\in \f (\r)$ such that

(3.9) \qquad  $\widehat {f} (\la)=1$ on $[-(b+\delta),b+\delta]$ and  supp $\widehat {f} \st (-(b+2\delta), b+2\delta)$.

\noindent Let $\phi\in L^{\infty}(\r,X)$.  Set

(3.10)\qquad $\phi_1= \phi*f$ and $\phi_2 =\phi -\phi*f$.

\noindent Then  (3.9)  and Lemma 2.3 (ii) give

(3.11)\qquad $sp(\phi_1) \st $ supp$\, \widehat{f} \st(-(b+2\delta), b+2\delta)$ and $sp(\phi_2) \cap [-b,b] =\emptyset$.

\proclaim{Theorem 3.5} Let  $(T(t))$ be a holomorphic $C_0$-semigroup  on $X$ with generator $A$ such that sup $_{t \ge 0} ||T(t)|| < \infty$ and let $\phi\in L^{\infty}(\r,X)$ with  Beurling spectrum  $sp(\phi)$ satisfying   $i \,sp(\phi) \cap \sigma (A)= \emptyset$.  Then there is $x_{\phi}\in X$ such that the equations (1.1), (1.2) with $x=x_{\phi}$ have a bounded uniformly continuous mild solution on $\r$.
\endproclaim

\demo{Proof} Let $\phi_1, \phi_2$ be given by (3.10).

 Case $\phi_1$: Let $K_1= sp (\phi_1)$, $i\,K_2= \sigma (A) \cap i\r$ and $\eta= d(K_1,K_2)$.
  Then $K_1, K_2 \st (-c, c)$, where $c$ is defined  above (3.9). By Lemma 2.4 (iv) there exist $W$, $\va$  and a  bounded sequence of trigonometric polynomials $(\pi_n)$ such that $\pi_n (t)\to  \phi_1(t)$ as $n\to \infty$ uniformly on each bounded interval of $\r$ and the Fourier exponents of
$(\pi_n)$ belong to $ W$,  where  $W$ is
an open set $W\supset K_1$  with $d(K_2,W)\ge  d(K_2,\text{\,supp\,}\va)\ge \eta/2$,  $\va \in \h(\r)$
               and $\va = 1$ on $W$.
  Set $ Q (\la) = \va (\la) \, R(i\la,A)$ on $\r\setminus K_2$ and $0$ on $K_2$. Since  supp $\va \cap K_2 =\emptyset$,   $ Q$ is infinitely differentiable on $\r$ and has compact support, with  $\widehat {Q}  \in   C_b (\r, L(X)) \cap L^1 (\r, L(X))$ by partial integration, so with Proposition 2.5
           the $Q$  is the Fourier transform of  a bounded continuous function $F\in L^1 (\r,L(X))$. We claim that
 $u_{\phi_1} (t)= (F * \phi_1) (t)$ is a mild solution on $\r$ of the equation $u' (t)= Au(t) +\phi_1 (t)$, $u(0)= (F * \phi_1) (0)$.
  Using (3.8)and $sp(\pi_n) \st W$, one can show that $u_n (t) :=(F*\pi_n) (t)$ is a classical solution on  $\r$ of the equation  $u'(t)= Au (t)+ \pi_n(t)$, $u (0) :=(F*\pi_n) (0)$ and therefore is a mild solution.   So  by Lemma 1.1, $u_n (t)= T(t-t_0)(F*\pi_n  (t_0))+ \int_{t_0}^t T(t-s)\pi_n (s)\, ds$ for each  $t \ge t_0\in \r$, $u_n(0)= (F*\pi_n)(0)$, $n\in \N$. Passing to the limit when $n\to \infty$, we get  with Lemma 2.4 (i), (ii) and the dominated convergence theorem

$(F* \phi_1) (t)= T(t-t_0) (F*\phi_1)(t_0)+ \int_{t_0}^t T(t-s) \phi_1 (s)\, ds$ for each $t \ge t_0\in \r$.

\noindent Now Lemma 1.1 shows that $u_{\phi_1} (t)$ is a mild solution on $\r$ of the equation (1.1) with $\phi=\phi_1$ and $x= (F* \phi_1) (0)$.

 Case $\phi_2$: By (3.11)  case $\phi_2$, one has   $sp(\phi_2) \cap
  [-b,b] =  \emptyset$. It follows $P\phi_2\in C_{ub} (\r,X)$ by [9, Corollary 4.4 valid also for $X$-valued functions]. By Lemma 2.3,  $sp (\phi_2)\st  sp (P\phi_2) \st sp(\phi_2)\cup \{0\}$.  With $h\in \f (\r)$ such that $\widehat {h} (\la)=1$ on $[-a/2,a/2]$ and  supp $\widehat {h} \st (-a,a)$, one has $sp ((P\phi_2)*h)\st \{0\}$, $(P\phi_2)*h= x_0\in X$ [3, Theorem 4.2.2], so  by Lemma 2.3 (ii)  the
   spectrum  of $\psi : = P\phi_2-(P\phi_2)*h$ satisfies $sp (\psi) \st sp(\phi_2)$.
   By Lemma 2.4 (iii) there exists a sequence of trigonometric polynomials $(\pi_n)$ such that $\pi_n (t)\to  \psi(t)$ as $n\to \infty$ uniformly on each bounded interval of $\r$ and the Fourier exponents of
$(\pi_n)$ belong to $\r\setminus [-b,b]$. Using (3.8), one can show that $u_n (t) := (G*\pi_n) (t)$ is a classical  and then a mild solution  on  $\r$ of the equation  $u'(t)= Au (t)+ \pi_n(t)$, $u (0) :=G*\pi_n (0)$, where $G=G_3$ defined by (3.7) case $m=3$.   So  by Lemma 1.1 (1.4), $u_n (t)= T(t-t_0)(G*\pi_n  (t_0))+ \int_{t_0}^t T(t-s)\pi_n (s)\, ds$ for each  $t \ge t_0\in \r$, $u_n(0)= G*\pi_n(0)$, $n\in \N$. Passing to the limit when $n\to \infty$, we get  with Lemma 2.4 (i), (ii), $G\in L^1(\r,L(X))$ by Lemma 3.4(iii)

$(G* \psi) (t)= T(t-t_0) (G*\psi)(t_0)+ \int_{t_0}^t T(t-s) \psi (s)\, ds$ for each $t \ge t_0\in \r$.

\noindent By Proposition 2.1 (i) we conclude that $G*\psi,  G*\phi_2 \in C_{ub}(\r,X)$. By (2.4),   $(G*\psi)  \in  C^1(\r,X) $. By  Lemma 1.1  and [1, Proposition 3.1.15 valid also for $J=\r$]), $G*\psi$ is a classical solution on $\r$ of the equation   $w'(t)= Aw(t)+\psi(t) $, $w(0)=G*\psi (0)$.
  So,  $G*\psi(t)\in D(A)$ for each $t\in \r$. Since $G*\psi= G*P\phi_2- G*x_0$ and  $G*x_0=0$ by Lemma 3.4 (v), it follows  $G*P\phi_2(t)\in D(A)$ for each $t\in \r$, with Proposition 2.1  then

$(G*\psi)' (t)= G*\phi_2(t)= A(G*(P\phi_2-x_0) (t)) +P\phi_2 (t)-x_0 = A (P (G*\phi_2) (t)) + $

 $A (G*(P\phi_2) (0))+ (P\phi_2) (t) -x_0= -x_0 + A(G*(P\phi_2)) (0))+ A P(G*\phi_2)(t)+ P\phi_2 (t)$,

\noindent with $P(G*\phi_2)(t)\in D(A)$ for each $t\in \r$.  It follows that $G*\phi_2 (0)=  -x_0 + A(G*P\phi_2) (0))$. Hence  $v=G*\phi_2 $ is a mild solution of $v'(t)= Av(t)+ \phi_2 (t)$ on $\r$, $v(0)=G*\phi_2 (0) $ by  (1.3).

Finally,  $u_{\phi}:= F*\phi_1+G*\phi_2$ is a mild solution  of (1.1)    on $\r$ with
$x=u_{\phi}(0)= F*\phi_1(0)+G*\phi_2(0)=x_{\phi}$ bounded and uniformly continuous by the above.
  \P
\enddemo

\proclaim{Theorem 3.6} Let $A$,  $\phi$ be as in Theorem 3.5 and let $\phi_1, \phi_2$ be as in (3.10). There exist $ F, G\in L^1 (\r,L(X))$ so that

(i) For each $x\in X$, $ u (t)=T(t)x+F*\phi_1- F*\phi_1(0)+G*\phi_2-G*\phi_2(0)$ is the unique mild solution on $\r_+$ of (1.1). Moreover, $u\in C_{ub}(\r_+,X)$.

(ii) In addition, if $J=\r$, $(T(t))$ is a  bounded  $C_0$-group (sup $_{t\in \r} ||T(t)|| < \infty$), then for each $x\in X$ equations (1.1),(1.2) have a unique mild solution on $\r$ given by  $ u (t)=T(t)x+F*\phi_1- F*\phi_1(0)+G*\phi_2-G*\phi_2(0)$. Moreover,  $u\in C_{ub} (\r,X)$.
\endproclaim

\demo{Proof}
(i) By [1, Proposition 3.1.16] the solutions of $v'=Av$  on $\r_+$ are given by $v(t)= T(t)x$,
(i) follows by Theorem 3.5.

 (ii)
   With Theorem 3.5 only   $u' = Au$, $u(0) = x$ , has to be considered.

   Uniqueness: (1.4) gives  $u(0)=T(0-(-n))u(-n) = T(n)u(-n), u(-n) =
   T(-n)u(0)$, $u(t) = T(t-(-n))T(-n)u(0) = T(t)u(0)$, $t > -n$, $n \in \N $.

   Existence:  $u(t) : = T(t)x $, $t \in\r$, gives  $u(t_0) = T(t_0)x$ or $x =T(-t_0)
   u(t_0)$,  so $u(t) = T(t)T(-t_0)u(t_0) = T(t-t_0)u(t_0)$, with Lemma 1.1
   $u$ is a mild solution   on $\r$.    \P
\enddemo

\head{\S 4.  Existence  of generalized almost periodic solutions of equation (1.1) in the non-resonance case}\endhead

For $\A\st L^1_{loc}(\r,X)$ we define mean classes $\m\A$  by ([5, p. 120, Section 3 ])

(4.1)\qquad$\m\A:=\{f\in L^1_{loc}(\r,X):   M_h f\in  \A ,\,\, h >0\}$, where

 \qquad \qquad \,\,\, $M_h f (t)=(1/h)\int_0^h f(t+s)\, ds$.

\noindent Usually  $\A  \st  \m\A\st \m^2\A\st \cdots$  with  the $\st$ in general strict (see [6, Proposition 2.2, Example 2.3]).

  We denote by $\n$ any class of functions having the following properties:

(4.2)  $\n$ linear $\st L ^1_{loc}(J,\Bbb{X})\st X^{\r}$.

(4.3) $(\phi_n)\st \n \cap  C_{ub}$ and $\phi_n \to \psi $ uniformly on $\r$ implies $\psi\in \n$.

(4.4) $\phi\in \n$, $a\in \r$ implies $\phi_a \in \n$.

 (4.5) $ B\circ\phi \in \n$ for each $B\in L(\Bbb{X})$, $\phi\in \n\cap C_{ub}$.

\proclaim{Lemma 4.1}  If   $\n$ satisfies (4.2)-(4.5) and  $\phi\in L^{\infty}(\r,X)\cap \m\n$, $F\in L^1(\r,L(X))$ respectively $L^1 (\r,\cc)$ then $F*\phi\in \n \cap C_{ub}(\r,X)$.
\endproclaim

\demo{Proof} By Proposition 2.1 (i), $F*\phi$ exists and $\in C_{ub}$.
To $F$ there is a sequence of $L(X)$-valued step-functions $H_n= \sum_{j=1}^{m_n} B_j \chi_{[\al_j,\beta_j)} $  with
$||F - H_n|| _{L^1}  \to 0$,  so  $||F*\phi  -  H_n*\phi||_{\infty} \le ||\phi||_{\infty} ||F - H_n|| _{L^1} \to 0 $ as $n\to \infty$. With (4.2), (4.3)
it is enough to show  $(B \chi_I)*\phi \in \n$, $I = [\al,\beta)$ for each $B\in L(X)$. With [1, Proposition 1.6]
one has $(B \chi_I)*\phi = B(\chi_I *\phi)$, with (4.5)  we have to show
$\psi_{\al,\beta} : = \int_{\al}^{\beta} \phi( \cdot - s) \,ds = \int_{-\beta}^{-\al}
\phi( \cdot + s)\,ds = \int _0 ^{- \al} \phi( \cdot + s) \,ds  -  \int_{-\beta} ^0 \phi( \cdot + s) \,ds \in \n$.
Now  $\phi \in\m\n $  gives  $\int_0^h \phi( \cdot + s) \,ds  \in \n$  if  $h>0$, (4.4) for $a = -h$
then  $\int_0^h \phi(\cdot - h + s)\, ds =  - \int_0^{-h} \phi( \cdot + s)\, ds  \i\n$, which
gives $\psi_{\al,\beta} \in\n$.
 The proof of $F\in L^1 (\r,\cc)$ is the same. \P
\enddemo

Examples of $\n$ satisfying (4.2)-(4.5) are the  spaces of almost periodic  functions $AP=AP (\r,X)$, almost automorphic functions $AA(\r,X)$, Bochner almost automorphic functions $BAA(\r,X)$,
          bounded Levitan almost periodic functions $LAP_b(\r,X)$ [7, p. 430], linear subspaces with (4.2) of  bounded recurrent functions  of $REC_b(\r,X)=RC$ of [7, p. 427], Eberlein almost  periodic functions $EAP(\r,X)$ and so on.

\proclaim{Theorem 4.2} Let $A$,  $\phi$ be as in Theorem 3.5 and $\phi\in \m\n$
with $\n$ satisfying (4.1), (4.2)-(4.5). Then there is $x_{\phi}\in X$ such that the equation (1.1), (1.2) with $x=x_{\phi}$ has a   bounded uniformly continuous mild solution on $\r$ which belongs to $\n $.
\endproclaim
\demo{Proof} By Theorem 3.5,$u_{\phi}:= F*\phi_1+G*\phi_2$ is a mild solution  of (1.1)    on $\r$ with
$x=u_{\phi}(0)= F*\phi_1(0)+G*\phi_2(0)=x_{\phi}$,  bounded and uniformly continuous. Since $\phi\in \m\n \cap L^{\infty}$, $\phi_1= \phi*f \in \n \cap C_{ub}$   by Lemma 4.1; since $\n \cap C_{ub}  \st \m\n$   by [6, Proposition 2.2, p.1011],
          $\phi_1  \in \m\n$. Again with Lemma 4.1 one gets
 $F*\phi_1\in  \n \cap C_{ub}$. Similarly, $\phi_2 =\phi- \phi_1 \in \m\n \cap L^{\infty}$ and so $G*\phi_2\in \n \cap C_{ub}$ by Lemma 4.1. \P
\enddemo

We should remark for example that if $\phi$ is  bounded Stepanoff $S^p$-almost periodic, the $u_{\phi}\in AP$  [5,  (3.8), p. 134]. Also, if $\phi$ is a Veech almost automorphic function [19], then $u_{\phi}$ is a uniformly continuous Bochner almost automorphic function [21, p. 66], [7, (3.3)].

\proclaim{Example 4.3} $X = Y^n$, $Y$ complex Banach space, $A =$ complex $n \times n$ matrix,
$u$, $\phi$  $X$-valued in (1.1), $\phi\in L^{\infty}(\r,X)$. Then, if $sp (\phi)$ contains no purely imaginary eigen-value of $A$ and $\phi \in \m\n$, $\n$ with (4.1)-(4.5), then (1.1) has a   mild solution on $\r$ which belongs to $ \n \cap C_{ub}$. This extends a well known result of Favard [13, p. 98-99].
\endproclaim

 Another example would be a result on the almost periodicity of all
solutions of the inhomogeneous wave equation in the non-resonance
case [20, p. 179, 181 Th\'{e}or\`{e}me III.2.1], [1, Proposition 7.1.1], here one has a $C_0$-group, all
solutions of the homogeneous equation are almost periodic.

\indent School of Math. Sci., P.O. Box No. 28M, Monash University,
 Vic. 3800.

\indent E-mail "bolis.basit\@sci.monash.edu.au".

\indent Math. Seminar der  Univ. Kiel, Ludewig-Meyn-Str., 24098
Kiel, Deutschland.

\indent E-mail "guenzler\@math.uni-kiel.de".
\enddocument

\ref\no13\by  Levitan B. M. and  Zhikov V. V.,:  Almost
Periodic Functions and
 Differential Equations, Cambridge Univ. Press, 1982
\endref

As the author realized in his "Added to the proofs" in [15, p. 1272], many of his results are
erroneous without additional assumptions. Since  his additional assumption (2.1) (below) is far too restrictive and he does not
point out what statements are affected
and how they should be modified, this "Added to the proofs" is not really
helpful. However, in  his corrigendum  [16] submitted  to JDE, the author stated clearly that all his claims in [15] on the non-uniformly continuity are omitted and the uniform continuity of all functions in [15] is assumed.

In this note we show that there is a  misconception in [15] which we believe resulted from some other hidden  misconceptions in previous works. We hope that this note will help to understand better the concept  of reduced spectra.
  In section 2, we give some comments (A)-(D) showing that the main classes of almost periodic and asymptotically almost periodic functions do not satisfy (2.1).  In section 3 we give two Examples supporting our arguments in the comments. Example 3.1 was first used by Levitan [17, pp. 144, 212-213]. Here we give a new treatment of example 3.1 to show that $\phi $ is a Bochner almost automorphic function without using the result of Levitan [17, p. 144] that $\phi$ is a Levitan almost periodic function.
 Example 3.2 is very important because it shows  in particular that there is $\psi \in AP$  with derivative $\psi'$  continuous and bounded but $\psi'$ is not even recurrent or Poisson stable (see the  definitions  (3.1), (3.2) in \S 3). This seems to be a missing example in the literature. In section 4, we prove a spectral criteria  for the  bounded mild  solutions  of the inhomogeneous abstract Cauchy  problem

 (*) $\frac{d u(t)}{dt}= A u(t) +  \phi (t) $, $u(0)=x\in \Bbb{X}$, $t\in  {J}$,

where $A$ is a closed linear operator   on $\Bbb{X}$ and $\phi\in L^{\infty} ( {J}, \Bbb{X})$, $ {J} \in\{\r_+,\r\}$.

\noindent This criteria is useful particularly in the case when  $\phi$ is not necessary  uniformly continuous but $u$ is uniformly continuous.

 \head{\bf\S 1. Notation,  Definitions and Preliminaries}\endhead

In the following $J, \tilde {J} \in\{\r_+, \r\}$, where $\r_+ =[0,\infty) $, $\N_0=\N\cup \{0\}$, $\Bbb{X}$ is a complex  Banach space, $L(\Bbb{X})$ is the Banach space of linear bounded operators $B:\Bbb{X}\to \Bbb{X} $; the elements of all the other spaces ($\st \Bbb {X}^J$, not $\Bbb {Y}$ of \S 2  (B),(C)) considered   are functions
       $\phi : J \to\Bbb{ X}$ (not equivalence classes), the $= , + ,$ scalar
      multiplication are pointwise on $J$ (not a.e.), correspondingly $L^{\infty}(J,\Bbb{X})$
     has the norm  $||\phi||_{\infty}  : =$ sup $\{||f(t)|| : t  \in J \}$ (not  essential supremum),
  $BC(J,\Bbb{X})$ is the Banach space of $\Bbb{X}$-valued bounded continuous functions on $J$, $BUC(J,\Bbb{X})$ is the Banach space of $\Bbb{X}$-valued bounded uniformly continuous functions on $J$, $AP(\Bbb{X}):=AP(\r,\Bbb{X}) $ is the Banach space of $\Bbb{X}$-valued almost periodic functions on $\r$, $BAA(\Bbb{X}) $ (respectively $VAA(\Bbb{X})$) is the Banach space of $\Bbb{X}$-valued Bochner  (respectively Veech) almost automorphic functions on $\r$ [12, Definition 2], [19,  Definition 1.2.1, p. 722], [9, p. 430], $C_0(J,\Bbb{X})$ is the Banach space of $\Bbb{X}$-valued  continuous functions on $J$  vanishing at infinity, all with sup-norm $||\cdot||_{\infty}$.
The Schwartz space of rapidly decreasing $C^{\infty}$ functions on $\r$ will be denoted by $\f (\r)$. The Fourier transform of $f\in L^1 (\r,\cc)$ will be denoted by $\hat{f} (\la)=\int_{-\infty}^{\infty} e^{-i\la t} f(t)\, dt$,  $\gamma_{\la}$ respectively $\frak{g}$ will denote the functions $\gamma_{\la} (t)=e^{i\la t}$ respectively $\frak{g}(t)= e^{it^2}$, $t, \la\in \r$.

For the convenience of the reader we collect some further definitions, assumptions and needed earlier results,
          for $ \n  \st  \Bbb{X}^J  $.

$Invariant$ :  $\phi_a  \in   \n$ if  $\phi  \in  \n$, $a  \in  \r$ with  $translate$ $\phi_a (t):=\phi(t+a)$.

$Positive$-$invariant$: translate  $\phi_a  \in   \n$  if  $\phi  \in
\n$ and  $0\le a <\infty$.

$BUC-invariant$ : $\phi  \in BUC (\r,\Bbb{X})$ and $\phi|\, J  \in  \n$
imply $\phi_a|\, J \in \n$   for all $a  \in\r$.

$Uniformly\,\, closed$ : $\phi_n  \in  \n$, $n \in  \N$, and $\phi_n \to
\phi$ uniformly on $J$ implies $\phi \in \n$.

 (1.1)\qquad $\m\n (J,\Bbb{X})= \{ \psi\in L^1_{loc}(J,\Bbb{X}): M_h \psi \in \n, \, h >0\}$,

 \qquad \qquad $\m_*\n (J,\Bbb{X})= \{ \psi\in L^1_{loc}(J,\Bbb{X}): (M_h \Psi)|\, J \in \n, \, 0\not =h \in \r\}$,

 \qquad \qquad \qquad where $\Psi:=\psi$ on $J$, $\Psi=0$ on $\r\setminus J$,

(1.2) \qquad $M_h \Psi(t)=(1/h)\int_{0}^h  \Psi(t+s)\, ds$, $0\not= h\in\r$, $M_h \psi=M_h\Psi|\, J$, $h >0$.

(1.3) $(\Delta)$ :  $\phi  \in   L^1_{loc}(J,X)$, $\Delta_h \phi  \in   \n$
for  $h>0$  implies  $\phi - M_k \phi \in\n$
                                 for  $k > 0$.

(1.4)  $\n$ linear $\st L ^1_{loc}(J,\Bbb{X})$, $\n$ uniformly closed, $\n \, BUC$-invariant  ([10, (3.1)]).

(1.5) (i) $\gamma_{\la}\phi \in \n$ for each $\gamma_{\la} (t)=e^{i\,\la t}$, $\phi\in \n$ ([4, $(l_2)$, p. 60 ]),

\qquad \,\, (ii) $\n$ contains all constant functions ([4, $(l_3)$, p. 60]),

\qquad \, \, (iii) $ B\circ\phi \in \n$ for each $B\in L(\Bbb{X})$, $\phi\in \n$ ([4, $(l_5)$, p. 60]).

\noindent  The $spectrum$ of an $\phi\in L^{\infty} (J,\Bbb{X})$ $with$ $respect$ to a class $\n \st L^ 1 _{loc} (\tilde {J},\Bbb{X})$  with $\tilde {J} \st J$  is defined by ([2, Definition 4.1.2, p. 20], [4, p.118], [10, Definition 3.1], [13, Definition 3], $\hat{f}=$ Fourier transform)

  (1.6) \qquad  $ sp_{\n}  (\phi)  : = sp_{\n}  (\Phi)= \{\lambda \in \r : f \in L^1 (\r,\cc), \Phi*f  |\tilde {J}  \in \n$ implies

 \qquad \qquad \qquad \qquad \qquad \qquad \qquad \qquad \qquad \qquad \qquad \qquad \qquad \qquad  $  \hat {f} (\la) =0  \}$.

\noindent Here $\Phi =\phi$ on $\tilde {J}$, $\Phi=0$ on $\r\setminus \tilde {J}$. $sp_{\n}(\phi)$ is always closed in $\r$. For $J=\r \not = \tilde {J}$ the $sp_{\n} (\phi)$ of (1.6) coincides with the definitions in [2], [4], [9], [13] by (1.7).

(1.7)  If $\n\st  L ^1_{loc}(J,\Bbb{X})$ satisfies (1.4) and $\phi \in L^{\infty}(J,\Bbb{X})$, then

\qquad   $sp_{\n}(\phi) = sp_{\n}(\psi)$ for any   $\psi \in L^{\infty}(\r,\Bbb{X})$ with  $\psi=\phi$
on $J$.

\noindent ([4, Lemma 1.1 (C)]).

(1.8) If $\n\st  L ^1_{loc}(J,\Bbb{X})$ satisfies (1.4),  $\phi \in L^{\infty}(J,\Bbb{X})$,  $f\in L^1(\r,\cc)$ and $\Phi$ is defined as in (1.6),  then \,\,
 $sp_{\n}(\Phi*f) \st sp_{\n}(\phi) \cap $ supp $\hat{f}$.

  ([4, Corollary 2.3 (C)]]).

(1.9) If $\n \st BUC (J,\Bbb{X})$ satisfies (1.4), then $\n\st \m\n$

 \noindent ([6, Proposition 2.2]).

(1.10) If $\n \st L_{loc}^1 (J,\Bbb{X})$ satisfies $(\Delta) $ and $r\n \st \n$ for real $r >0$, then

 \qquad $\m\n\st \n+\n'$,

 \noindent where  $\n' (J,\Bbb {X})=\{\psi: \psi=\phi' \text {a.e.\, with\,} \phi\in\n (J,\Bbb {X}) \cap W^{1,1}_{loc} (J,\Bbb {X})\}$.

 \noindent ([9, Proposition 1.1, p. 38]).

(1.11) If $\n$ is convex uniformly closed $\st BUC(J,\Bbb{X})$, then $\n$ satisfies $(\Delta)$.

\noindent ([6, p. 1012, Proposition 3.1], [8, Theorem 2.4, p. 428]).

(1.12)\qquad $(\Delta)$ for linear positive  invariant $\n$ implies $\n\st \m\n$.

(1.13) If  $\n \in\{AP(\Bbb{X}),   VAA (\Bbb{X}), C_0 (J,\Bbb{X})\}$, then $\n$ satisfies  $\n\st\m\n$ and $(\Delta) $.

\noindent See (1.11), (1.12),  [8, Proposition 3.5 (ii), p. 431].

(1.14) Let $\n$ be a positive-invariant  closed linear subspace $\st BC (J,\Bbb{X})$. If $\phi\in \n$ and $\phi'$ is uniformly continuous, then $\phi'\in \n$.

 \noindent ([2, Proposition 1.4.1], [6, Proposition 2.9]).

(1.15) Let $A: D(A)\to  \Bbb{X} $ be a closed linear  operator with linear domain $D(A)\st \Bbb{X}$. Let $I \st \r$ be an interval (bounded or unbounded) and let  $h: I\to \Bbb{X} $ be a Bochner integrable function with $h(t)\in D(A)$ for each $t\in I$. Suppose that $A\circ h$ is Bochner integrable. Then $\int_I h(s)\, ds \in D(A)$ and $A \int_I h(s)\, ds=\int_I Ah(s)\, ds$.

 \noindent  [1, Proposition 1.1.7, p. 11].

 \proclaim {Proposition 1.1}  For any  $\n  \st  \Bbb{X}^{\tilde {J}}$,  $\phi  \in  L^{\infty} (J,\Bbb{X}) $, if $\phi \in \m\n$ then  $sp_{\n} (\phi)  = \emptyset$.
\endproclaim

 \demo{Proof} For any $\lambda \in\r$, define  $h = \pi/|\lambda|$ if $\lambda \not = 0$, else
$h = 1$; then the step function $f=(1/h) \chi_{(-h,0)}  \in  L^1(\r,\r)$ and with
$\Phi : = 0$ outside $J$,$\Phi = \phi$ on $J$,  one has $f*\Phi |J = M_h \phi  \in\n$ ,
with  $\hat{f}(\lambda) \not = 0$, so  $\lambda  \not \in sp_{\n} (\Phi)$.  It follows $sp_{\n} (\Phi)=\emptyset $ and so $sp_{\n} (\phi)=\emptyset $ by (1.6). $\square$
\enddemo
\proclaim {Corollary 1.2} If  $\phi  \in\n  \st\m\n$   and $\phi  \in  L^{\infty}(J,\Bbb{X})$, then  $sp_{\n} (\phi)
                 =  \emptyset $. This is false without $\n\st \m\n$ by $\phi = \frak{g}$ and $\n=\frak{g}\cdot AP (\r,\cc)$.
\endproclaim

\proclaim {Proposition 1.3} Let   $\phi  \in  L^{\infty} (J,\Bbb{X}) $ and let $\Phi $ be defined as in (1.6). Let   $\n  \st L^1_{loc} (J,\Bbb{X})$ satisfy (1.4)  and $\la_0 \in sp_{\n} (\phi) $.

(i) There is $h_0 >0$ such that $\la_0 \in sp_{\n} (M_{h_0}\phi) $.

(ii) For any $f\in L^ 1(\r,\cc)$ with $\hat{f} (\la_0)\not = 0$, one has  $\la_0 \in sp_{\n} (\Phi*f) $.
\endproclaim

 \demo{Proof} (i) With   (1.7), this is a special case of [10,(3.11), (4.1)].

(ii) Assuming $\lambda_0 \not \in  sp_{\n} (\Phi*f)$, there is $h \in L^1 (\r,\cc)$  with $\hat{h} (\la_0)\not = 0$ and $\Phi*(f*h) |J=  (\Phi*f)*h |J \in\n$, with
     $\widehat {(f*h)}(\lambda_0)\not = 0$; this implies  $\lambda_0  \not \in sp_{\n} (\phi)$,
      against the assumption. $\square$
\enddemo

In the following we identify $L^1(I,\cc)$ respectively $\n \st L^{\infty} (I,\Bbb{X})$ with the subspace $\{f\in L^1(\r,\cc): f(t)=0, t \in \r\setminus I\}$ respectively $\{\phi \in L^{\infty}(\r,\Bbb{X}): \phi|\,I \in \n, \phi=0, t \in \r\setminus I\}$ . Here $I\in\{\r, \r_+, \r_-\}$, $\r_- = (-\infty, 0]$.

\proclaim {Proposition 1.4} Let  $\n $ be a linear uniformly closed subset $ of L^{\infty} (J,\Bbb{X})$. Then in (a), (b) the  conditions (i), (ii) are equivalent.

(a) (i) $\n\st \m\n$, \qquad  (ii) $\n * L^1 (\r_-,\cc)|\,J \st  \n$.

(b) (i)  $\n\st \m_*\n$, \qquad  (ii) $\n * L^1 (\r,\cc)|\,J \st  \n$.

(c) If $\n$ satisfies even (1.4), then (a)(i) is also equivalent with (b)(ii).
\endproclaim

\demo{Proof} (a)   $(i)\Rightarrow (ii)$: With $\Phi=\phi$ on $J$, $\Phi=0$ on $\r\setminus J$, we have

         (1.16)\qquad  $M_h \phi   =   (1/h)(\Phi *\chi_{(-h,0)})|\, J$,
          where   $\chi_{(a,b)} : = 1$ on  $(a,b)$,  $: = 0$ on  $\r\setminus J$.

   \noindent        As  $M_h \phi =(\Phi*s_h)|\, J\in \n$, $\phi\in\n $, $h >0 $, $s_h=(1/h) \chi_{(-h,0)}$, it follows $\Phi*\xi|\, J \in\n $  for all
               step functions  $\xi$  on $\r_-$; since these are dense in
               $L^1(\r_-,\cc)$ and $\n $ is  uniformly closed, (ii) follows.

\noindent  $(ii)\Rightarrow (i)$ follows by (1.16) and $s_h \in L^1(\r_-,\cc)$ for each $h >0$.

(b)The proof is similar to part (a).

(c)  By  (a)(i), (1.16) and  the property that $\n$ is $BUC$-invariant, one has

$(M_h \phi)_t|\, J =(\phi*s_h)_t|\, J= \phi*(s_h)_t|\, J\in \n$, $\phi\in\n $, $h >0 $, $t\in \r$

\noindent and so $\phi*\zeta|\, J \in\n $  for all
               step functions  $\zeta$  on $\r$; since these are dense in
               $L^1(\r,\cc)$ and $\n $ is uniformly  closed, (b)(ii) follows.

\noindent  $(b)(ii)\Rightarrow  (a)(i)$ follows by (1.16). $\square$
\enddemo
\proclaim {Proposition 1.5} Let  $\n $ be a linear uniformly closed subset of $ L^{\infty} (J,\Bbb{X})$.

(a) $\n\st \m\n$ implies   $\n \cap BUC(J,\Bbb{X})$ is positive invariant.

(b) $\n\st \m_*\n$  implies   $\n \cap BUC(J,\Bbb{X})$ is $BUC$-invariant.

(c) If $(\n * E)|\, J\st \n$ for some dense subset $E\st L^1(\r,\cc)$ implies   $\n \cap BUC(J,\Bbb{X})$ is $BUC$-invariant.
\endproclaim

\demo{Proof}
 Let $\phi\in L^{\infty} (J,\Bbb{X})$. Then

(1.17)\qquad $M_{h+k}\phi = \frac{h}{h+k}M_{h}\phi+ \frac{k}{h+k}(M_{k}\phi)_h$,\,\, $ h+k \not =0$, $ h\not =0$, $ k \not =0$.

\noindent (a) If $\phi \in \n$, $ h+k  > 0$ and $h > 0$, then $M_{h+k}\phi,\, M_{h}\phi\in \n$  since $\n\st \m\n$.  It follows  $(M_{k}\phi)_h\in \n$, $k + h>0$, $h > 0$  by linearity of $\n$ and (1.17).  As $\lim_{k\to 0} M_k\phi_h= \phi_h$ uniformly on $J$, $\phi\in (\n \cap BUC(J,\Bbb{X}))$, we get $\phi_h \in (\n \cap BUC(J,\Bbb{X}))$ for each $h >0$ and $\phi\in (\n \cap BUC(J,\Bbb{X}))$. This gives (a).

(b) This  follows as in (a) from (1.17),
     $0 \not= h  \in\r$    fixed, $k \to 0$.

(c) The assumptions imply $\n * L^1 (\r,\cc)|\,J \st  \n$. So, (c)
 follows by Proposition 1.4 (b) and part (b).
  $\square$
\enddemo
 An "inverse" of Proposition 1.5 is false,
 $\n\st \m\n$ does not imply

 $BUC$-invariance:

\proclaim{ Example 1.6} $\n =\frak{g}AP (\r,\Bbb{X})$  is linear, uniformly  closed, invariant, $\st BC (\r,\Bbb{X})$,  but  $\n\not\st \m\n$ ($\frak{g}\in \m C_0(\r,\cc)$, $(\frak{g}AP)\cap C_0 (\r,\cc)=\{0\}$). See also Example 4.4 (ii).
\endproclaim

\proclaim{ Example 1.7} $\n  : = \{ \phi  \in   C_0(\r,\r): \phi = 0$  on $\r_+ \}$ is
    linear , uniformly closed, positive invariant,  $\st BUC(\r,\r)$, with  $\n\st \m\n$,
    but $\n$ is not invariant.

\noindent Here $\n\not\st \m_*\n$, $\m_* \n$ strictly $\st \m\n$.
\endproclaim

 \head{\bf\S 2. Explicit comments  }\endhead

In the following we give some comments which had been written before the author kindly informed us about his corrigendum [16]. So, the comments (A), (B), (D) are concerned with the now to be omitted false results about the cases of non-uniformly continuous functions.

(A) The induced operator $\widetilde{\h}$ of  [15, Definition 2.6, respectively p. 1263 for $J=\r_+$] is well defined for  a class $\n \st BC(J,\Bbb{X})$ satisfying Condition  $F$ or  Condition  $F^+$ of [15, Definitions 2.3, 3.1] if and only if

\qquad (2.1) \qquad $\h(D(\h)) \cap \n) := \{f ' : f  \in  C^1(J,\Bbb{X}) \cap \n, f'\,$ bounded$\} \st  \n$.

\noindent This is possible if $\n =\{0\}$ or $\n = BC(J,\Bbb{X})$. We are not aware of other interesting classes satisfying Condition $F$ or $F^+$ and (2.1):

(B) If $\n\in \{ BUC(\r,\cc), AP (\r,\cc), BAA (\r,\cc)\}$ respectively $\n=  C_0 (J,\cc)$, then $\n$ satisfies [15, Definition 2.3] respectively [15, Definition 3.1], but does not satisfy (2.1)
 and  so the operator $\widetilde{\h}$  is not well defined by part (A):

Indeed,
 if $\n\in \{ BUC(\r,\cc), AP (\r,\cc), BAA (\r,\cc)\}$, $\Bbb{Y}= BC(\r,\cc)/\n $, $\Phi (t)=\int_0 ^t \psi(s)\, ds$ with $\psi$ of Example 3.2 below,   then $ \Phi \in \n\cap D(\h)$ but $\Phi'= \psi$ bounded $ \not\in \n$.

  If
 $\n =C_0 (J,\cc)$, $\Bbb{Y}= BC(J,\cc)/\n $, $\phi (t) = \int_0^1 \frak{g} (t+s)\, ds$, then $ \phi \in \n\cap D(\h)$ but $\phi'$ is bounded  with $\phi' (\sqrt{2\pi n})\not \to 0   $ as $n\to \infty$. Here $\frak{g} (t)= e^{it^2}$.

If  $\n\not = \{0\}$ satisfies  Condition $F$ of [15], then $AP(\r ,\Bbb{X})\st \n$ and if $\n$ satisfies  Condition $F^+$, then $C_0 (J,\Bbb{X})\st \n$. So,  either $\n$ does not satisfy (2.1) or one can show that

(2.2)\qquad $ (\m AP \cap  BC(\r,\Bbb{X})) \st \n$ respectively $(\m C_0 (J,\Bbb{X})\cap BC(J,\Bbb{X})) \st \n$.

\noindent Indeed, by  the statements (1.13),(1.11), (1.10),  we have $\m\n =\n+\n'$ for $\n=AP$ or $C_0$,  implying (2.2).

Since $\widetilde{\h}$ is not well defined for the known  classes of interest, Definition  2.8 of [15] makes no sense if $\n \not =\{0\}$. So with the approach of  [15] the reduced spectra  $sp_{\n}(\phi)$, $sp_{\n}^+(\phi)$ of [15, Definitions 2.8, 3.3] are not  defined. Therefore the inclusions [15, (4.3), (4.9)] have no meaning. However, see Theorem 4.3 below.

(C) If $\n \st BUC(J,\Bbb{X})$  satisfies (1.4),
  $\Bbb{Y}= BUC(J,\Bbb{X})/\n$ and $D(\h)=\{\psi\in BUC(J,\Bbb{X}): \psi'\in BUC(J,\Bbb{X})\}$, then the induced operator $\widetilde{\h}$ is well defined by (1.14).
  In this case [15, Definitions 2.8, 3.3]
for $\phi\in  BUC(\r,\Bbb{X})$ give the Carleman spectrum of the class $\tilde{\phi}=$
the Beurling  spectrum  (see [10, (3.4)]), coinciding with the usual  $sp_{\n} (\phi)$ defined in (1.5) by a result of
          Chill and Fasangova (see [13], [10, Theorem 3.10]) and [1, Proposition 4.8.4, p. 321].

With Examples 3.1, 3.2, the assumptions (1.4) are optimal and not " very strict" as  the author claimed in [15, p. 1252, 1263], $BUC(J,\Bbb{X})$, $BC(J,\Bbb{X})$ satisfy (1.5). Contrary to the footnotes 1), 3) of [15, p. 1250, 1263], $sp_{\n} (\phi) =\emptyset$ is admissible in [9]. In [16] the author realized that his approach works only for classes $\n \st BUC(J,\Bbb{X})$ satisfying Condition $F$ or $F^+$.  By Proposition 1.4 (c) and [11, Proposition 2.1], it follows that if  $\n$ satisfies  Condition $F$ or Condition $F^+$, then $\n$ is $BUC$-invariant. This means that if $\n$ satisfies  Condition $F$, then $\n$
satisfies   (1.4) and (1.5) and if $\n$ satisfies Condition $F^+$, then $\n$ satisfies (1.4), (1.5)(i), (iii).

 In  (D) we mean results  of [15] but we use $sp_{\n}(\phi)$ as defined  in  (1.6), (1.7).

(D)   Corollary 2.20 as stated is not correct, already for $\n = AP(\r,\cc)$:
    This $\n$ satisfies Condition $F$ of [15, Definition 2.3], but by Examples 3.1, 3.2 below there exists $\psi  \in  BC(\r,\cc)$
     with $sp_{AP} (\psi) = \emptyset$, but $\psi  \not \in AP(\r,\cc)$.
 Theorem 2.25, Corollaries   2.27 and 2.29  are not correct for $sp_{AP}$
       also by Example 3.1 below ; for $sp_{\n}$ with $\n=\{0\}$ this is open.  Corollary  2.31 is not correct by Example 3.2 below.

 Theorem 3.6, Corollaries 3.8, 3.9 are not true:  If $\phi =\int _0^1 \frak{g} (\cdot +s)\,ds \in C_0 (\r_+,\cc)$  (see part (B)), then $\phi' = (\frak{g} (\cdot+1)- \frak{g})\in BC(\r_+,\cc)$ and $sp_{C_0} (\phi')=\emptyset$ by Proposition 1.1 but $\phi'\not\in C_0 (\r_+,\cc)$.

\head{\bf\S 3. Two examples}\endhead

Here  follow the two examples used to refute many of the  results of [15]. Here we give a direct proof for example 3.1. Example 3.2 is instructive for various conclusions concerning wide classes of almost periodicity (for example almost periodic (ap), almost automorphic (aa), Levitan almost periodic (L-ap), recurrent and Poisson stable functions). For the benefit  of the reader we give the relevant definitions.

(3.1) By a $recurrent$ function $\phi$  we mean  $ \phi\in REC(\r,X) : = \{\phi  \in  C(\r,\Bbb{X})  : E(\phi,1/n,n)$
relatively\,\, dense in $\r$ for each $n \in  \N \}$, with

 $E(\phi,\e,n) : = \{ \tau \in \r :
||\phi(t+\tau)-\phi(t)||\le \e\,\,$ for all  $\,\,|t|\le n\}$.

$E(\phi,1/n,n)$  is $relatively\,\, dense$ means there is a compact
set $K \st \r$ such that $K + E (\phi,1/n,n)= \r$ (see [15, Definition 2, p. 80],  [8, p. 427]).

(3.2) An $\phi\in C(\r,\Bbb{X})$ is $Poisson\,\, stable$ if it has at least one sequence $(t_m)\st \r$ with $t_m\to\infty$  such that $\phi_{t_m}\to \phi$ locally uniformly in $\r$
(see [15, Definition 1, p. 80]).

\proclaim{Example 3.1} The function $\phi= \sin \frac {1}{p}\in BC(\r,\r)$ with $p (t) = 2+\cos t+\cos \sqrt{2}t$ is Stepanoff almost periodic $S^1$-$AP\st \m AP$ and  $sp_{AP}(\phi)=\emptyset$ but $\phi\not\in AP= AP (\r,\cc)$.
See [4, p. 119, (1.5), p. 118 above  (1.2), (3.6), (3.9)] for the definitions. This $\phi$ is also Bochner almost automorphic (B-aa) [12] and so Veech almost automorphic  (V-aa)  [19] and L-ap [8, p. 430, (3.4)] (see also [4, p. 119] and references therein).
 \endproclaim

\demo{Proof} First we show that  $\phi\in S^1$-$AP$. Set
$\phi_n (t): = sin \frac{1}{2+\text{ max\,\, }\{\cos \,t, -1+\frac{1}{n}\}+\cos \sqrt{2}t}$. Then $\phi_n\in AP$ for each $n\in \N$ and $\phi_n (t)= \phi (t)$  if $ \text{ max\,\, }\{\cos \,t, -1+\frac{1}{n}\}= \cos\, t$. It follows $\int _0^{2\pi} |\phi_n (t+s)-\phi (t+s)|\, ds\le  2\mu (E^t_n)$, where $\mu$ is the Lebesgue measure on $\r$ and $E^t_r =\{\tau\in [t,t+2\pi]:  \text{ max\,\, }\{\cos \,\tau, -1+\frac{1}{r}\}= -1+\frac{1}{r}, \, r \ge 1 \}$.  Then $\mu (E^t_n)  =\mu (E^0_n)  = \mu ([\pi - \delta_n,\pi + \delta_n])$  with $\cos \delta_n
           = 1 - 1/n$, $t \in\r$,  with $\delta_n \to 0$ as $n \to \infty$.
    It follows $\lim_{n\to\infty} \int _0^1 |\phi_n (t+s)-\phi (t+s)|\, ds=0$ uniformly  in  $t\in \r$ and implies  $\phi\in S^1$-$AP$ ( see [4, p. 132]).
  So $M_h \phi (\cdot)= (1/h)\int_0^h \phi (\cdot +s)\, ds \in AP$  for each $h >0$ by [4, (3.8)]. By Proposition 1.1,  one gets
 $sp_{AP}(\phi)=\emptyset$.

  Now, we show that
                  $\phi$ is not uniformly  continuous. Indeed, since range of $p$ is  $R(p)= (0,4]$, for each $n\in \N$, by Kronecker's approximation
               theorem [14, p. 436, (d)] and continuity,  there is  $t_n >0$ such that
$p (t_n)= \frac{1}{n\pi}$.
Choose $t'_n$ nearest point to $t_n $ with $p (t'_n)= \frac{1}{(n-\frac{1}{2})\pi}$. We have $|t_n-t'_n|\le  \mu (E^{t_n}_{(n-\frac{1}{2})\pi})\to 0$ as $n\to\infty$.
 Since $|\phi(t'_n)-\phi(t_n)|=1$, we get  $\phi $ is not uniformly  continuous.
    It follows  $\phi\not \in AP $.

  Finally, we show that   $\phi$ is B-aa. Indeed,  since $(\cos\,t , \cos\sqrt{2} t)$ is almost periodic,
   for $(t_{n'})\st \r$  there are $\al$, $\beta \in \r$ and  a subsequence
$(t_{n})$ such that

(3.4) \qquad $\cos(t+t_n) \to \cos(t+\al)$, $\cos (\sqrt {2}(t+t_n)) \to \cos (\sqrt {2}(t+\beta))$,

\qquad\qquad  $p(t+t_n) \to (2 +\cos(t+\al)+ \cos (\sqrt {2}(t+\beta)))=: q(t)$, uniformly in $t\in\r$.

\noindent Since $q$  is entire, $
C :=\{s\in \r: q(s)=0\}$ is at most  countable.
  So, there is a (diagonal)  subsequence $(s_n)$  and $\psi:\r\to [-1,1]$ with $\phi (t+s_n)=  \sin \frac {1}{p (t+s_n)} \to \psi(t)$ pointwise for each $t\in\r$. Now, (3.4) implies $q(t-s_n)\to p (t)$,  $p(t+s_m-s_n)\to q(t-s_n)$ and then $p(t+ s_m-s_n)\to p (t)$ as $(n,m)\to \infty$ for each $t\in \r$.
This yields $\phi (t+ s_m-s_n)\to \phi (t)$ as $(n,m)\to \infty$ pointwise in $t\in \r$; $m\to\infty$ and the definition of $\psi$ give therefore  $ \psi(t-s_n)\to \phi (t)$. By  definition 2 of [12], $\phi$ is B-aa.

    See also [17, pp. 212-213] for another proof that $\phi \in S^1$-$AP$ but $\phi\not \in AP$; and [8, Example 3.3] that $\phi$ ia B-aa.
    $\square$
\enddemo

 \proclaim{Example 3.2}  There is $\psi  \in  BC(\r,\r)$ which is not ap or B-aa  or V-aa or
    recurrent or uniformly continuous (not even Poisson stable  (see (3.1), (3.2) respectively Example 3.1),  also $\Delta_1 \psi (\cdot): =\psi(\cdot+1)- \psi(\cdot)$ and so $\psi$  are not Stepanoff $S^1$-almost periodic),  but  $P\psi(t): =\int_0^t  \psi(s)\,ds$  is almost periodic and so $sp_{AP}(\psi)= sp_{BAA}(\psi)=\emptyset$.
\endproclaim

\demo{Proof}   Take  $\psi =  \sum_{n=1}^{\infty} h_n$,\qquad $h_n$ periodic  with period $2^{n+1}$,

     $ h_n (t)=  0  $, $\,\,\, t \in [-2^n, 2^n  - 1]$,\qquad
       $h_n(t)= \sin (2^n \pi t )$, $\,\,\, t\in [2^n  -1,2^n]=: I_n$.

\noindent   One has
 $\text{\, supp\,} h_{n}= I_n + 2^{n+1}\Bbb{Z}$ and  for each $n\not =m$,
 $\text{\, supp\,} h_{n}\cap \text{\, supp\,} h_{m} =\emptyset$: The right endpoints of the translations of  $I_n$ are all
         even, and if  $n = m+k$, $k \in \N_0$, $2^n + 2^{n+1}u = 2^m +2^{m+1}v$
        implies  $ 2^k (1+2u) = (1+2v)$  and then  $k=0$, $u=v$.
  It follows $\psi\in BC(\r,\r)$ and  with $I= [-2,0]$ for each $\tau \ge 2$, $r\in \N$

(3.5) \qquad $\text{\, sup \,}_{t\in  I} |\psi (t+\tau)-\psi (t)|= \text{\, sup \,}_{t\in  I} |\psi(t+\tau)|  \ge \text{\, sup \,}_{t\in  I} |h_r(t+\tau)|$.

   \noindent Since  $\int_{I_n} h_n (t)\, dt=0$, $Ph_n $ is periodic with period $2^{n+1}$ and $||Ph_n||_{\infty}\le 2^{-n}$.  It follows $P\psi \in AP (\r,\r)$.  This implies $M_h \psi  \in  AP$ for  $h > 0 $, and so  $sp_{BAA} (\psi) \st
                   sp_{AP} (\psi) =  \emptyset$  by Prop. 1.1.

 With $ \delta = 2^{-n}$ one has

           $\int_{I_n} |\Delta_1 \psi (t+\delta) - \Delta_1 \psi (t)|\,dt  \ge
           \int_{I_n}|\psi(t+\delta)-\psi(t)|dt  - \int_{I_n} |\psi(t+1+\delta)-\psi(t+1)| dt
          = $

          $ \int_{I_n} |h_n(t+\delta)-h_n(t)|dt - \int_{I_n}|\psi(t+1+\delta)|dt \ge
         2/\pi  - 2^{-n}$ ,  $>  0.1$  for all  $n  \in  \N$.

         \noindent It follows  $\Delta_1 \psi$ and so $\psi$ are not uniformly continuous even in the $S^p$-norm (see [4, p. 132] for the definition).  Hence $\psi, \Delta_1 \psi$ are not almost periodic  and not $S^1$-almost periodic.

 Since to each even $ n  \in  \N$  there exist unique  $m\in \N$,  $k\in \N_0 =\N \cup  \{0\}$ such that $n= (1+2k)2^m$, we get

  (3.6) \qquad           $n  =   2^m +  k  2^{m+1} $.

\noindent  We show that
   for  each
     $\tau  \ge 2$ there is $r\in \N$ with
      \,\, $\text{\, sup \,}_{t\in  I} |h_r (t+\tau)| =1$.
 Indeed,  let $\tau \in 2\N+ y$ for some $y\in  [0,2]$.  Then by
  (3.6), since $2n+y = 2(n+1)+y'$ with $y'=y-2 $

  $\tau= 2^m+k 2^{m+1}+y$ for unique $m\in \N$, $k\in \N_0$  and $y\in  [0,\frac{3}{2}]$ or

  $\tau= 2^{m'}+k' 2^{m'+1}+y'$ for unique $m'\in \N$, $k'\in \N_0$  and $y'\in  [-\frac{1}{2},0]$.

\noindent With $t = -y - 2^{-m-1}$ respectively  $t = -y' - 2 + 2^{-m'-1}$ we get

  $  \text{\, sup \,}_{t\in  I} |h_m(t+2^m + k 2^{m+1}+y)| = 1 $ for each $y\in [0,\frac{3}{2} ]$,

  $  \text{\, sup \,}_{t\in  I} |h_{m'}(t + k' 2^{m'+1}+2^{m'}+y')| = 1 $ for each $y'\in [-\frac{1}{2}, 0]$.

 \noindent  By (3.5),  it follows $\text{\, sup \,}_{t\in  I} |\psi(t+\tau)-\psi(t)|\ge 1$ for all  $\tau \ge 2$.
 Since B-aa  and V-aa functions are always recurrent, we conclude
$\psi$ is not B-aa or V-aa or recurrent or Poisson stable by the definitions  (3.1), (3.2).
$\square$
\enddemo

\head {\bf \S 4. Reduced  spectrum of solutions of evolution equations}\endhead

In this section we study  the reduced  spectrum with respect to  a class $\n\st L^1_{loc} (\tilde {J},\Bbb{X}) $ of bounded solutions of evolution equations

(4.1) \qquad $\frac{d u(t)}{dt}= A u(t) +  \phi (t) $, $u(0)=x\in \Bbb{X}$, $t\in J$,

 \qquad\qquad where $A$ is a closed linear operator   on $\Bbb{X}$ and $\phi\in L^{\infty} (J, \Bbb{X})$.

 \noindent The  half-line  (Laplace) spectrum denoted by  $sp_L(\psi)$ for $\psi\in L^{\infty} (\r_+,\Bbb{X})$ is introduced in [1, p. 275]. If  $\n \st L^1_{loc} (\r_+,\Bbb(X))$ satisfies (1.4), then
 $sp_{\n} (\psi) \st sp_{w} (\psi)\st sp_{L} (\psi)$, by [10, (3.12), (3.14)]. Here $ sp_{w} (\psi)$ is the weak half-line  (Laplace) spectrum [1, Definition 4.9.1, p. 324]. The reduced spectrum and the half-line spectrum
 of solutions of (4.1) when $u, \phi\in BUC(J,\Bbb{X})$ have been investigated by many authors  see for example  [3], [1, sections 4, 5] and lists of references there. In this section we  prove inclusions (4.5), (4.6) for (4.1) which are known for the
 half-line  spectrum  of solutions of (4.1) in the case $u, \phi\in BUC(\r_+,\Bbb{X})$, see [1, Proposition 5.6.7 b), p. 380].

  \proclaim{Definition 4.1 (see [1, p. 120, 121,  380])}
  A function $u\in  C(J,\Bbb{X})$ is called a mild solution of (4.1) if $\int_0^t u(s)\, ds \in D(A)$, $x\in\Bbb{X}$ and $u(t)-x = A \int_0^t u(s)\, ds+ \int_0^t \phi (s)\, ds$, $t\in J$.
  \endproclaim

\proclaim{Proposition 4.2} Let $u\in  BC(\r,\Bbb{X})$ be a mild solution of (4.1), $J=\r$, with $\phi\in L^{\infty} (\r, \Bbb{X})$. If  $g\in \f (\r)$, $n\in\N_0$, then

   \qquad \qquad $v =u*(g^{(n)})= (u*g)^{(n)}$ is a classical  and so a mild solution of

  (4.2)\qquad  $\frac{d v(t)}{dt}= A v(t) +  \phi*g^{(n)} (t) $, $v (0)=u*g^{(n)}(0)$, $t\in \r$.

  \noindent Moreover, if $\n \st L^1_{loc} (\tilde {J},\Bbb{X})$ satisfies (1.4), (1.5) and  supp $\hat {g}\cap  sp_{\n} (\phi)=\emptyset$, then

   (4.3)\qquad $i\, sp_{\n} (v)\st \sigma (A)\cap i\,\r $.
\endproclaim
\demo{Proof} We have $P(u*g)(\tau)= ((Pu)* g)(\tau) =\int_{\r} (\int_0^{\tau} u(t-s)\, dt) g(s)\, ds$, by Fubini's Theorem [1, Theorem 1.1.9, p. 12].  Since $D(A)$ is linear,  Definition 4.1 gives  $ g(s)\int_0^{\tau} u(t-s)\, dt\in D(A)$ for each $\tau, s\in\r$. Since $ A \int_0^{\tau} u(t-s)\, dt$ is with Definition 4.1  $s$-continuous and  $O(|s|)$, one gets $ (A \int_0^{\tau} u(t-s)\, dt) g(s)$ is Bochner integrable on $\r$ for each $\tau\in\r$. By  (1.15), it follows

(4.4)\qquad  $ A P(u*g)(\tau)= ((A Pu)*g) (\tau) - ((A Pu)*g)(0)$.

 \noindent This and Definition 4.1 for $u$ proves that  $v= u*g$ is a mild solution of (4.2) for the case $n=0$.  Since $g\in\f (\r)$, $v=u*g$, $\phi*g$ are infinitely differentiable, with $(u*g)'= u*(g')$. It follows   that $v$ is a classical solution  by an extension  of [1, Proposition 3.1.15, p. 120] to $\r$. This means $v*g(t)\in D(A)$, $t\in\r$. Since $g^{(n)}\in \f(\r)$ for each $n\in \N_0$, (4.2) follows.

If supp $\hat {g}\cap  sp_{\n} (\phi)=\emptyset$, then  $sp_{\n} (\phi*g)=\emptyset$, by (1.8). Since
$\phi*g \in BUC(\r,\Bbb{X})$, we conclude $\phi*g|\, \tilde {J} \in \n$  by [4, Corollary 2.3(A)].  By the above, we have $v= u*g\in BUC(\r,\Bbb{X})$ is a classical solution of (4.2) satisfying $v(0), v'(0)\in D(A)$. The assumptions imply that $ \n\cap BUC(\r_+,\Bbb{X})$  is  a $\Lambda$-class satisfying [3, $(l_1) $-$(l_5)$].
 The inclusion (4.3) follows by [3, Theorem 3.3 (valid also for $A$ only closed)].
 $\square$
\enddemo

  \proclaim{Theorem 4.3} Let $\n \st L^1 _{loc}(\tilde {J},\Bbb{X})$ satisfy (1.4), (1.5)  and  let  $\phi \in L^{\infty}(J,\Bbb{X})$.   If $u\in BC(J,\Bbb{X})$ is a mild solution of (4.1), then

  (4.5)\qquad $i\,sp_{\n} (u) \st ((\sigma (A)\cap i\,\r)\cup i sp_{\n} (\phi))$.

\noindent If moreover   $\phi\in \m\n$, then

(4.6)\qquad  $i\,sp_{\n} (u) \st \sigma (A)\cap i\,\r$.

\endproclaim
\demo{Proof} First we prove the case $J=\r$.
 Let  $u\in BC(\r,\Bbb{X})$ be a mild solution of equation (4.1).
  Let $\rho (A)$ be the resolvent set of $A$ and let $ i\la_0 \in U= (\rho (A)\cap  i\r) \cap (i\r \setminus i sp_{\n} (\phi))$. Since $U$ is an open set, there is $1 >  \delta >0$ and $f\in \f(\r)$ such that $i (\la_0-\delta, \la_0 +\delta)\st U$ and supp $\hat{f} \st (\la_0-\delta, \la_0+\delta)$. It follows supp $\hat{f} \cap sp_{\n} (\phi) =\emptyset$. By (4.3) and  (1.8), we have $i sp_{\n} (u*f)\st \sigma (A)\cap i\r$,   $i sp_{\n} (u*f)\st i (\la_0-\delta, \la_0+\delta)$.
  As $i (\la_0-\delta, \la_0+\delta) \cap  (\sigma (A)\cap i\r) =\emptyset$, we get $i sp_{\n} (u*f)=\emptyset$. Since $u*f\in BUC(\r,\Bbb{X}) $, we conclude $u*f|\,J \in \n $ by [4, Corollary 2.3 (A)], and so $i\la_0\not \in sp_{\n } (u)$.  This proves (4.5).
 If $\phi \in\m\n$, then  $sp_{\n} (\phi) =\emptyset$  by  Proposition 1.1. This and (4.5) give (4.6).

The case $J=\r_+$.
Let  $u\in BC(\r_+,\Bbb{X})$ be a mild solution of equation (4.1) and let $k, h >0$. With  [1, Proposition 3.1.15 p. 120] for $M_h u$, one can  show that  $v_{k,h}= M_k M_h u$ is a classical solution of

(4.7) \qquad  $\frac{d v(t)}{dt}= A v(t) + M_k M_h\phi (t) $, $v(0)= M_k M_h v (0)$, $t\in \r_+$ with

(4.8) \qquad   $ v_{k,h}' (0)\in D(A)$.

\noindent Moreover,  $v_{k,h}, \psi_{k,h} = M_k M_h\phi \in BUC(\r_+,\Bbb{X})$. Let  $V_{k,h}, \Psi_{k,h} \in BUC(\r,\Bbb{X})$ be $C^1$-extensions of $v_{k,h}, \psi_{k,h}$ satisfying

\qquad  $\frac{d V(t)}{dt}= A V(t) + \Psi(t) $, $V_{k,h}(0)= M_k M_h v (0)$, $t\in \r$.

\noindent This is possible by [3, Lemma 3.2 (i)] and (4.8).  Applying the case $J=\r$, we get

(4.9) \qquad  $i\,sp_{\n} (V_{k,h}) \st ((\sigma (A)\cap i\,\r)\cup i sp_{\n} (\Psi_{k,h}))$.

 \noindent Using (1.6), [10, Lemma 4.2, (4.1), (3.11)] by (1.8) and (4.9), with $\Phi$  as defined in (1.6) we conclude

   $i\, sp_{\n} (u)= \cup_{k >0, h >0}i\, sp_{\n} (M_k M_h u)= \cup_{k >0, h >0} i\, sp_{\n} (V_{k,h})\st $

 \qquad   $\cup_{k >0, h >0} ((\sigma (A) \cap\, i\r)\cup \, i\, sp_{\n} (\Psi_{h,k}))=$
  $ (\sigma (A) \cap\, i\r)\cup $

  $(\cup_{k >0, h >0}(i\, sp_{\n} (M_h M_k \Phi))= (\sigma (A) \cap\, i\r)\cup i\,sp_{\n} (\Phi)= (\sigma (A) \cap\, i\r)\cup i\,sp_{\n} (\phi)$.

 \noindent This proves (4.5). Finally, (4.6) follows  from (4.5) by the same argument as in the case $J=\r$.
$\square$
\enddemo
In the following example we consider the case $\Bbb{X}=\cc$ and $A: \cc\to \cc$  defined by $Ac=ic$. We have $\sigma(A)=\{i\}$. Example 4.4 (ii) shows that the condition $\phi\in \m\n$ in  Theorem 4.3 (ii) can not be replaced by $\phi\in \n$. Also, it supports the suspicion  that  (4.5) might be trivial  without the condition $\n\st \m\n$. Example 4.4 (ii)   shows also that  (1.4), (1.5)  do not imply $\n\st \m\n$. Example 4.4 (iii), (iv) show
that (4.6) can be valid  though the conditions (1.5) (i), (ii) are  not satisfied. One can optimize conditions (1.5) (i), (ii) using [10, Theorem 4.3].

\proclaim{Example 4.4} All solutions of the equation $y'(t)=iy (t)+\frak{g}(t)$ are given by

$y(t)= e^{it} (c+ \int_{0}^{t} e^{-is}\frak {g}(s)\, ds)=: \gamma_1 (t)( c + y_1(t))$, where $c\in \cc$.

(i) If $\n = \frak{g}\cdot AP(\r,\cc)$, then $sp_{\n}(\frak{g})=\r$ and so (4.5) is trivially satisfied.

(ii) If $\n = \frak{g}\cdot AP(\r,\cc)\oplus AP(\r,\cc)$, then  $\n$ satisfies (1.4), (1.5) but  $sp_{0}(\frak{g})=sp_{\n}(y)=\r$ and so (4.6) is not satisfied.

(iii) If $\n = (\frak{g}\cdot AP(\r_+,\cc))\oplus (C_0 (\r_+,\cc)\oplus (\gamma_1\cdot \cc))$, then $\n\st \m\n$,  $sp_{\n}(\frak{g})=\emptyset$ and $sp_{\n}(y) = \emptyset$, $c\in \cc$.

(iv) If $\n = (\frak{g}\cdot AP(\r,\cc))\oplus (\gamma_1\cdot C (\overline{\r},\cc))$, then $sp_{\n}(\frak{g})=sp_{\n}(y)= \emptyset$, $c\in \cc$.

Here $C(\overline{\r},\cc)=\{\phi\in BUC(\r,\cc): \lim_{t\to\infty} \phi (t), \lim_{t\to -\infty} \phi (t) \text {\,\, exist}\} $.
\endproclaim

\demo{Proof} (i)  we have $\frak{g} \in \n =\frak{g}\cdot AP$ but $sp_{\n} (\frak{g}) = sp_{\n \cap BUC (\r,\cc)} (\frak{g})=sp_{0} (\frak{g})= \text{supp\,} \hat {\frak{g}}=\r$ (see  [4, (1.3)], [10, (3.4)], [10, Example 2.2]). It follows (4.5) is trivially satisfied.

 (ii)   $\n$ satisfies (1.4), (1.5), but we omit the proof that  $\n$ is uniformly closed.
  We have $\frak{g} \in \n$ but $sp_{\n} (\frak{g}) = sp_{\n \cap BUC (\r,\cc)} (\frak{g})=sp_{AP} (\frak{g})=\r$ since $\frak{g}* f \in AP(\r,\cc)$ implies  $\frak{g}* f \in AP(\r,\cc)\cap C_0 (\r,\cc)=\{0\}$, $f\in L^1 (\r,\cc)$,  and  $sp_{0} (\frak{g})=\r$.

  Now we show that $sp_{\n} (y)=\r$. We have  $sp_{\n} (y) = sp_{\n \cap BUC (\r,\cc)} (y)= sp_{AP} (y)$
   and
$sp_{AP}(\gamma_1 y_1)\st sp_{AP}   (y)  \cup  sp_{AP}(-c \gamma_1)$ [2, Theorem 4.1.4]. As
$sp_{AP}(-c \gamma_1) = sp_{AP} (\gamma_1) = \emptyset $ ([2, Theorem 4.1.4]),
  $sp_{AP} (\gamma_1 y_1 )$
 $ = 1 + sp_{AP} (y_1)$ ([2, Theorem 4.1.4]), $sp_{AP} (u') \st sp_{AP} (u)$
 ($u  \in  BUC (\r,\cc)$,
   $u'  \in   BC (\r,\cc)$), and  $sp_{AP} (y_1 ') = $  $ -1  + sp_{AP} (\frak{g})$,
 $sp_{AP} (\frak{g})=\r$
  shown already,
 it follows  $sp_{\n}(y)=\r$.

(iii)  We have $y,\,\frak{g}|\,\r_+ \in \n$ and $\n$
satisfies  $\n\st\m\n$. So, $sp_{\n}(\frak{g})=sp_{\n}(\frak{g})=\emptyset$ by Corollary 1.2.

(iv)   We have  $\frak{g}, y \in \n$ and  $\n\st\m\n$. The result follows by Corollary 1.2.    $\square$
\enddemo

\indent School of Math. Sci., P.O. Box No. 28M, Monash University,
 Vic. 3800.

\indent E-mail "bolis.basit\@sci.monash.edu.au".

\indent Math. Seminar der  Univ. Kiel, Ludewig-Meyn-Str., 24098
Kiel, Deutschland.

\indent E-mail "guenzler\@math.uni-kiel.de".

\enddocument

Thank you for your e-mail of  today.

1- I corrected the remark after Th. 5.2. But I kept also ref. [21]= [Pruss]. I think your objection was due to the misprint of the page occurred which led you to think that [21]=[Ph+S]=[20].

2. I confirm that exponential dichotomy is equivalent to hyperbolic $T(t_0)$ for some $t_0 >0$. This also in [EN, One-Parameter .., p.305, Prop. 1.15].

Intuitively notice that $sp(T(t))$ can be guessed from $sp(T(1))$ for $t >0$ for example if $\lambda \not\in sp(T(1)$ and $|\lambda|=1$, then $\lambda \not\in sp(T(t)$. [Also notice that $T(n)= (T(1))^n$].

3. I agree with your objection  against the uniqueness Theorem of yesterday. The remedy is the following:

In Lemma 5.1 of NRE we should add:

4. Then $v$ is the unique almost periodic  classical solution  of  (instead of 'Then $v$ is a  classical solution  of ')

5. In Theorem 5.2: At the end we should add:

Moreover $u= G*\phi$ does not depend on the choice of $H$ satisfying the assumptions of Lemma 4.3.

In the proof of Theorem 5.2, we add that $u_n:= G*\Pi_n $ are  unique  classical   and so mild solutions from $AP(\r,X)$ of (5.3) on $\r $ with (instead '$u_n := G*\Pi_n$ are   classical and so mild solutions  of (5.3) on $\r $ with '). Further we should use this to conclude that $u= G*\phi$ does not depend on the choice of $H$.

\proclaim{Corollary} Let  $ \n\in \{AP, AA, REC_b\}$, $A$, $\phi$ as in Theorem 5.2 and $\phi\in L^{\infty}\cap \m\n$. Then there exists $ G \in  L^1(\r,L(X))$ so that  the  $u _{\phi}: = G*\phi$ is the unique mild  solution on $\r$  of (5.3) which belongs to $\n \cap BUC(\r,X)$.
\endproclaim

6- I will change Lemma 5.1 and Th.5.2 only if you approve 4., 5.

Best wishes, Bolis.

\proclaim {Proposition 1.4} Let  $\n \st L^{\infty} (J,\Bbb{X})$ be linear and uniformly closed. Then in (a), (b) the  conditions (i), (ii), (iii) are equivalent.

(a) (i) $\n * L^1 (\r_-,\cc)|\,J \st  \n$, \,(ii) $\n\st \m\n$, \,    (iii)  $\n \cap BUC(J,\Bbb{X})$ is positive invariant.

(b) (i) $\n * L^1 (\r,\cc)|\, J \st  \n$, \, (ii) $\n\st \m_*\n$, \, \,  (iii)  $\n \cap BUC(J,\Bbb{X})$ is  $BUC$-invariant.
\endproclaim
\demo{Proof} (a) The implication  $(i)\Rightarrow (ii)$ follows from

         (1.15)\qquad  $M_h \phi   =   (1/h)(\phi *\chi_{(-h,0)})$ ,
          where   $\chi_{(a,b)} : = 1$ on  $(a,b)$,  $: = 0$ else in  $\r$.

   \noindent        $(ii)\Rightarrow (iii)$: Let $\phi\in L^{\infty} (J,\Bbb{X})$. Then

(1.16)\qquad $M_{h+k}\phi = \frac{h}{h+k}M_{h}\phi+ \frac{k}{h+k}(M_{k}\phi)_h$,\,\, $ h+k \not =0$, $ h\not =0$, $ k \not =0$.

\noindent If $\phi \in \n$, $ h+k  > 0$ and $h > 0$, then $M_{h+k}\phi,\, M_{h}\phi\in \n$ by (ii).  It follows  $(M_{k}\phi)_h\in \n$, $k + h>0$, $h > 0$  by linearity of $\n$ and (1.16).  As $\lim_{k\to 0} M_k\phi_h= \phi_h$, $\phi\in (\n \cap BUC(J,\Bbb{X}))$, we get $\phi_h \in (\n \cap BUC(J,\Bbb{X}))$ for each $h >0$ and $\phi\in (\n \cap BUC(J,\Bbb{X}))$.

 \noindent  $(iii)\Rightarrow (ii)$ follows by (1.11), (1.12).

 \noindent   $(ii)\Rightarrow (i)$:  By (1.15),   $\phi*\chi_{(-k,-h)} = \phi*\chi_{(-k,0)} -
                \phi*\chi_{(-h,0)} \in\n$  if  $0< h< k$, and so $\phi*S \in\n $  for all
               step functions  $S$ vanishing on $\r_+$; since these are dense in
               $L^1(\r_-,\cc)$ and $\n $ is uniformly  closed, (i) follows.

(b) The proof is similar to part (a).
  $\square$
\enddemo

Comment on the paper titled 'A spectral theory of continuous functions and the Loomis-Arendt-Batty-Vu theory on the asymptotic behavior of solutions of evolution equations' by  Nguyen Van Minh

\proclaim {Proposition 1.4} Let  $\n \st L^{\infty} (\r,\Bbb{X})$ satisfy (1.4). Then in (a) and (b) the  conditions (i), (ii), (iii) are equivalent.

 (i) $\n\st \m\n$, \,  (ii) $\n * L^1 (\r,\cc) \st  \n$, \,  (iii)  $R(\la, \Cal{D})\n \st \n$, $\la\in \cc$,  Re $ \la \not =0$.
\endproclaim

\demo{Proof}  By (i) and  the property that $\n$ is $BUC$-invariant, we have  $(M_h \phi)_t =(\phi*s_h)_t= \phi*(s_h)_t\in \n$, $\phi\in\n $, $h >0 $, $t\in \r$, where $s_h= (1/h)\chi_{(-h,0)}$. Since  $ {span }\{(s_h)_t: h >0 , t\in \r\}$ is dense in $L^1$ and $\n$ is uniformly closed, one gets  $\n * L^1 (\r,\cc) \st  \n$. This proves $(i)\Rightarrow (ii)$. Obviously,
$(ii)\Rightarrow (i)$ and so $(i)\Leftrightarrow (ii)$.

 With $\check {f}_{\la} (s) =f_{\la} (-s)$, $s\in \r$, we have

\

 $R(\la,A)\phi (t)= \cases{ \int_{0}^{\infty}e^{-\la s} \phi (s+t)\,ds=\phi*\check{f}_{\la}(t),
\text{\,\, if\,\,\,\,} \text{Re\,\,}\la > 0}\\ { -\int_{-\infty}^ 0 e^{-\la s} \phi (s+t)\,ds=\phi*\check{f}_{\la}(s),\text {\,\,\, if \,} \text{Re\,\,}\la <
 0}\endcases$,\,\,\, where

\

$f_{\la}(t) = \cases{ e^{-\la t},
\text{\,\, if\,} t \ge 0}
\\ { 0,\text {\qquad \,\,\, if \,} t<
 0}\endcases$, \,\, $ \text {Re\,} \la >0$, \,\,
 $f_{\la}(t) = \cases{0, \text {\,\, if \,} t>0
 }\\{- e^{-\la t},
\text{\, if\,\,}  t \le 0}\endcases$, \,\,\,$ \text {Re\,} \la  < 0$.

By  the property that $\n$ is $BUC$-invariant, we have $\phi*(f_{\la})_t= (\phi*f_{\la})_t\in \n$, Re $\la \not =0$, $\phi\in \n$.

\noindent Since $ {span }\{(f_{\la})_t: \text{ Re\,\,}\la \not =0, t\in \r \}$ is dense in $L^1$ and $\n$ is uniformly closed, one gets  $\n * L^1 (\r,\cc) \st  \n$. This proves $(iii)\Rightarrow (ii)$. Obviously, $(ii)\Rightarrow (iii)$ and so $(ii)\Leftrightarrow (iii)$. $\square$
\enddemo

\proclaim{Proposition A. J. Pryde} $ E= span \,\{f_{\la}: \text{ Re\,\,}\la \not =0\}$ is a dense  subspace of  $L^1 (\r,\cc)$.
\endproclaim
\demo{Proof} We have  $(L^1 (J,\cc))^{*} = L^{\infty}(J,\cc)$. If $E$ is not dense in $L^1 (\r,\cc)$, then by Hahn-Banach theorem there is
$0\not = \phi\in L^{\infty}(J,\cc)$ such that

$\Cal  {C}\phi(\la)= \int_0^{\infty} e^{-\la t}\phi (t)\, dt=0$, Re $\, \la >0$,

$\Cal  {C}\phi(\la)= -\int_0^{\infty} e^{\la t}\phi (-t)\, dt=0$, Re $\, \la  < 0$.

This means that the Carleman transform $\Cal  {C}\phi$ is zero on $\cc\setminus i\r$ and implies $sp_{C} (\phi)=\emptyset$ and so $\phi =0$. $\square$
\enddemo

Dear Hans,

I will try to answer  your questions:

1- You are right, $y\not\in AAP(\r,\cc)$ but $y\in AAP(I,\cc)$  with $I\in \{\r_+, \r_-\}$.

2-I received your e-mail that you can resolve this. But also, see Analysis Paper 122, Example 2.2. There you will find a beautiful  proof that $sp_{0}(\frak {g})=\r$ and  $sp_{L}(\frak {g})=\emptyset$ using $sp_{0}(\frak {g})=sp_{0}(\frak {g}_t)= sp_{0}(\frak{g}(t^2)\gamma_{2t}\frak {g})=2t+ sp_{0}(\frak {g})$ for each $t\in \r$ and $sp_{0}(\frak {g})\not = \emptyset$. Similarly for $sp_{L}(\frak {g})$.

3- For all $c\in\cc$, we have $0\not =y\in BUC(\r,\cc)$.  For all $f\in L^1$, we have $y*f \in BUC(\r,\cc)$ and so $y*f\not\in \frak{g}\cdot AP$. Also, $sp_0 (y)=\r$ because $sp_0 (\phi) \st sp_0 (P\phi)$ for each $\phi, P\phi \in BC(\r,\cc)$ and $sp_0 (\phi_1+ \phi_2) \st sp_0 (\phi_1) \cup sp_0 (\phi_2) $. This implies that for each $\la$ there $f\in L^1$ such that $\hat {f}(\la)\not =0$, $y*f\not =0$ implying $sp_{\n} (y)$, where $\n= \frak{g}\cdot AP$.

4- $J\in\{\r_+,\r\}$ but I withdraw this statement replace $J$ by $I\in \{\r_+, \r_-\}$. See also the new version of Ex. 4.4.

5- 8 see the new Ex. 4.4.

\

\proclaim {Proposition 1.3} Let  $\n  \st L^1_{loc} (J,\Bbb{X})$ satisfy (1.7). Let   $\phi  \in  L^{\infty} (J,\Bbb{X}) $, $\Phi =\phi$ on $J$, $\Phi=0$ $\r\setminus J$ and let   $\la_0 \in sp_{\n} (\phi) $.

(i) There is $h_0 >0$ such that $\la_0 \in sp_{\n} (M_{h_0}\phi) $.

(ii) For any $f\in L^ 1(\r,\cc)$ with $\hat{f} (\la_0)\not = 0$, one has  $\la_0 \in sp_{\n} (\Phi*f) $.
\endproclaim

 \demo{Proof} (i) Take   $h_0 = \pi/|\lambda|$ if $\lambda_0 \not = 0$, else
$h_0 = 1$.  Then $g_0=(1/{h_0}) \chi_{(-h_0,0)}  \in  L^1(\r,\r)$ and $\hat{g_0}(\lambda_0) \not = 0$.  One has $\Phi*g_0 |J = M_{h_0} \phi \not \in\n$ by (1.5), (1.6). We show that $\la_0 \in sp_{\n} (M_{h_0}\phi) $. Indeed, assuming  $\la_0 \not \in sp_{\n} (\Phi*g)$, there is $g  \in  L^1(\r,\r)$  with supp  $\hat {g}= V $ compact and  $\hat{g} (\la_0)\not =0$ such that    $sp_{\n} (\Phi*g_0) \cap V =\emptyset$.    It follows $ sp_{\n}(\Phi *g_0* g)=\emptyset$, by (1.9) (ii). So,  $\Phi *g_0* g \in \n$, by [4, Corollary 2. 3(A)] though  $k=g_0*g$ has  $\hat {k} (\la_0)\not = 0$.   This is a contradiction which proves   $\la _0 \in sp_{\n} (\Phi*g_0)$.

(ii) If $f\in \f(\r)$ with $\hat{f} (\la_0)\not = 0$, then $\Phi*f |J \not \in\n$ by (1.5). Arguing as in part (i), we get $\la_0 \in sp_{\n} (\Phi*f) $.
\enddemo
Case $J=\r$:  If $u\in BC(\r,\Bbb{X})$ is a mild solution of equation (3.1) with $J=\r$, then
   $v= M_h u$ \, is a  classical  solution of

(3.5) \qquad $\frac{d v(t)}{dt}= A v(t) + M_h \phi (t) $, $v(0)= M_h u(0)$, $t\in \r$.

  \noindent  Indeed,  since $\phi\in L^{\infty}(\r,\Bbb{X})$,   $M_h \phi\in BUC(\r,\Bbb{X}) $. From Definition 3.1, we conclude that     $v(t)\in D(A)$ for each $t\in \r$ and $\frac{d v(t)}{dt}= (u(t+h)-u(t))\in BC(\r,\Bbb{X})$. So $v\in C^1 (\r,\Bbb{X})$ and $v$ satisfies (3.5). It follows $v\in BUC (\r,\Bbb{X})$ and $v$ is a classical solution of (3.5) by [1, Proposition 3.1.16, p. 120].  Arguing in the same way as  in the proofs of [3, Theorem 3.3, Corollary 3.4] we conclude

   $i\,sp_{\n} (v) =i\, sp_{\n} (M_h u)\st ((\sigma (A) \cap\, i\r)\cup sp_{\n} (M_h \phi) )$.
   By [10, (4.1), (3.11)], it follows  $i\, sp_{\n} (u) = \cup_{h > 0}i\, sp_{\n} (M_h u) \st  ((\sigma (A) \cap\, i\r)) \cup \cup _{h > 0}i\, sp_{\n} (M_h \phi) ) = ((\sigma (A)\cap i\,\r)\cup i sp_{\n} (\phi))$. This proves (4.2).  If $\phi \in \n$, then $M_h \phi \in \n$ because $\n\st \m\n$. As  $M_h \phi \in BUC(\r,\Bbb{X})$, $sp_{\n} (M_h \phi)=\emptyset$ for each $h >0$. The inclusion (3.4) follows.

 Case $J=\r_+$. Let  $u\in BC(\r_+,\Bbb{X})$ be a mild solution of equation (3.1) with $J=\r_+$ and $f\in \f(\r)$.  let $U$, $\Phi$ denote respectively the functions $u$, $\phi$ extended by zero on $(-\infty,0)$. Consider the equation

 (3.6) \qquad $\frac{d v(t)}{dt}= A v(t) + (\Phi*f) (t) $, $v(0)= (U*f)(0)$ $t\in \r$.

\noindent Then $v= (U*f)$ is a mild solution of (3.6) by (1.7) and so it is a classical solution  by  [1, Proposition 3.1.5, p. 120]. The inclusions (4.2), (3.4) follow  by  the above and [4, Lemma 1.1, Corollary 2.3].

As  the reduced spectra  $sp_{\n}(f)$, $sp_{\n}^+(f)$ of [15, Definitions 2.8, 3.3] are not  defined for the main classes of interest, the inclusions [14, (4.3), (4,9)] are not   proved for $u\in BC(J,\Bbb{X})$. However, the inclusion (3.12) of [4] is sharper than [6, (4.9)]. Formula (4.3) of [6] can be deduced from the case when $u, f$ of [6, (4.1)] are uniformly continuous using the correct definition of reduced spectrum  pointed  in (C),

(2) \qquad $\frac{d v(t)}{dt}= A v(t) + M_h f (t) $, $t\in\r,  h > 0$ and

\qquad \qquad $sp_{\n} (g) = \cup_{h > 0} sp_{\n} M_h g$ for each $g\in L^{\infty} (\r,\Bbb{X})$.

\noindent Indeed, if $u$ is a mild solution of [6, (4,1)], then $v= M_h u$ is a mild solution of (2.1) and $v, M_h f\in BUC(\r,\Bbb(X))$.

one can show that  since $\m\n =\n+\n'$, where  $\n' (J,\Bbb {X})=\{g: g=f' \text {a.e.\, with\,} f\Bbb{X}\in\n (J,\Bbb {X}) \cap W_{loc^{1,1}} (J,\Bbb {X})\}$  by [2, [Proposition 3.21(i)], p. 678]

Dear Hans,

I answer some of your questions and this part will be deleted.

1- $\m AP \cap BC \st \n$ is valid because we assume that $\n$ satisfies (1). I think you mean $AP \st \m AP \cap BC $, which we  may mention without proof.

2.  Since $BAA_u \st BAA$, I think $BAA_u $ may cause confusion.

3. See the new version of Ex. 2.

4. I do not know, but I will think about that.

5. I think constructing  $\n$ with $F$ or $F^+$ and (1) should be the last thing to be concerned with. May be Minh should think about that.

Example 2 is now has the answer of 3. If you want we give reference [RJMP, Example 3.3] for $\phi$ of Example 1 to be BAA.

\

We recall that if $ H(t)=\sum_{r=1}^{\infty}f_r (t) $, where $F$, $(f_r)\st L^{\infty} (\r,\cc)$, then Beurling spectrum \qquad
$sp(F)\st \overline {\cup_{r=1}^{\infty} sp(f_r)}$.

Indeed, let  $t\in U:= \r\setminus \overline {\cup_{r=1}^{\infty} sp(f_r)}$.  Since $U$ is open  there exists $\delta >0$ such that $(t-\delta, t+\delta) \st U$. It follows $(t-\delta, t+\delta) \cap  sp(f_r)=\emptyset$ for each $r\in\N$. Choose $ \varphi\in \Cal {S} (\r)$ such that $\text{supp\,}\widehat {\varphi} \st (t-\delta, t+\delta)$, $\widehat {\varphi} (t)\not =0$. Then $sp (f_r*\varphi) =\emptyset$ by [6, 12 Lemma, p. 989] and so $f_r*\varphi =0$ for each $r\in\N$. This implies $F*\varphi =0$ and implies $t\not\in sp (F)$.

For a proof for  the case when the sum is finite  see [6, 11 Lemma (e), p. 988].


\proclaim {Proposition} $\phi= \sin \frac{1}{p}\not \in C_{ub}(\r,\r) $.
 \endproclaim

\demo{Proof} We have
Range of $p$ is  $R(p)= [4,0)$. For each $n\in \N$, by Kronecker's approximation
               theorem there is  $t_n >0$ such that
$p (t_n)= \frac{1}{n\pi}$.
Choose $t'_n$ nearest point to $t_n $ with $p (t'_n)= \frac{1}{(n-\frac{1}{2})\pi}$. We have $|t_n-t'_n|\le  \mu (E^{t_n}_{(n-\frac{1}{2})\pi})\to 0$ as $n\to\infty$.
 Since $|\phi(t'_n)-\phi(t_n)|=1$, we get  $\phi\not \in C_{ub}(\r,\r) $. $\square$
\enddemo






\proclaim   {Lemma} To each even $ n  \in  \N$  there exist $ m \in \N$ and  $k \in \N_0$  with
              $n  =   2^m  + k  2^{m+1}$,
(The  $m, k$ are then unique).
\endproclaim
\demo{Proof} Choose $m\in \N$ such that $n= (2k+1)2^m$ for some $k\in \N_0$. $\square$
\enddemo









   I- Proof that $\phi= \sin \frac{1}{p}\not \in C_{ub}(\r,\r) $ [4, p. 213]: Range of $p$ is  $R(p)= [4,0)$. For each $n\in \N$, there are $t_n,  t'_n >0$ such that $p (t_n)= \frac{1}{n\pi}$, $p (t'_n)= \frac{1}{(n-\frac{1}{2})\pi}$ and  $(t_n-t'_n)\to 0 $ as $n\to\infty$. We have $|\phi (t_n)-\phi (t'_n)|=1$ but $|t_n-t'_n|\to 0 $ as $n\to\infty$.

In  [4, p. 213], $t_n, t'_n$ are chosen to be $t_n =$ min $\{ t: p(t) =\frac{1}{n\pi}\}$, $t'_n =$ min $\{ t': p(t') =\frac{1}{(n-\frac{1}{2})\pi}\}$ and showed that $|t_n-t'_n|\to 0$ as $n\to\infty$.

One could also choose first $t_n$ with  $p (t_n)= \frac{1}{n\pi}$ then choose $t'_n$ nearest point to $t_n $ with $p (t'_n)= \frac{1}{(n-\frac{1}{2})\pi}$. Then $(t_n-t'_n)\to 0 $ as $n\to\infty$ because $p\in C_{ub}$.

\

II- I do not have a simple proof for

 3- The same idea can be used to show that $\text{\, sup \,}_{t\in  [-1,1]} |f(t+\tau)-f(t)| =1$ for each $\tau \ge 2 $.

  Also, I agree that $ h_k (\cdot+\tau)|\, [-1,1]\not  = h_k |\, [-1,1]$ is not enough unless supp  $ h_k (\cdot+\tau)|\, [-1,1]$ contains an interval of the form $[x,x+1/2] \st [-1,1]$.

   Your Lemma as stated is not correct

My argument is by cases and ultimately  we should use a correction of your Lemma:

Case $\tau \in [2,4]$: For $\tau \in [2,7/2]$ we calculate  supp $h_1 (\cdot+\tau)$ and show that $ h_1 (\cdot+\tau)|\, [-1,1]\not  =0$ and this is enough to show that
$\text {sup \,}_{t\in  [-1,1]} |f(t+\tau)-f(t)| =1$.

For $\tau \in [8/2, 4]$ we calculate  supp $h_2 (\cdot+\tau)$ and show that $ h_2 (\cdot+\tau)|\, [-1,1]\not  =0$ and this is enough to show that
$\text {sup \,}_{t\in  [-1,1]} |f(t+\tau)-f(t)| =1$.

Case $\tau \in [4,8]$ can be discussed similarly.

Case $\tau \in [2^k, 2^k+3/2]$ we calculate  supp $h_{k} (\cdot+\tau)$ and show that $ h_{k} (\cdot+\tau)|\, [-1,1]\not  =0$ and this is enough to show that
$\text {sup \,}_{t\in  [-1,1]} |f(t+\tau)-f(t)| =1$.

Case $\tau \in [2^k+ 3/2, 2^k+7/2]$, $k \ge 2$, we have $h_1 (t+\tau)= h_1 (t+ 2^k + x)=h_1 (t+x)$, where $x \in [3/2,7/2]$ we calculate  supp $h_1 (\cdot+x)$ and show that $ h_1 (\cdot+x)|\, [-1,1]\not  =0$ and this is enough to show that
$\text {sup \,}_{t\in  [-1,1]} |f(t+\tau)-f(t)| =1$.

Case $\tau \in [2^k+ 7/2, 2^k+11/2]$, $k \ge 3$, we have $h_2 (t+\tau)= h_2 (t+ 2^k + x)=h_1 (t+x)$, where $x \in [8/2,11/2]$ we calculate  supp $h_2 (\cdot+x)$ and show that $ h_1 (\cdot+x)|\, [-1,1]\not  =0$ and this is enough  to show that
$\text {sup \,}_{t\in  [-1,1]} |f(t+\tau)-f(t)| =1$. And so on.

I rewrite Definition 2.6 of [3]:

\proclaim{Definition 2.6} $D(\widetilde{\h})$ is defined to be the set of $\Bbb{Y}=BC(\r,X)/\n$ consisting of classes of functions that contain elements of $D(\h)=\{f\in BC(\r,X): f'\in BC(\r,X)\}$. The operator $\widetilde{\h}$ is defined by $\widetilde{\h}\tilde{f}=\tilde{f}'$ whenever $\tilde{f}\in D(\widetilde{\h})$ containing $f\in D(\h)$.
\endproclaim

Claim: $\tilde{f}'$ is not unique: Take the class $\tilde{f}=\n = AP(\r,\cc)$ which is the zero class of $\Bbb{Y}=BC/AP$. Take $\phi\in BC(R,\cc)$ with $\phi\in \m AP$ but $\phi\not \in AP$. Then $f=M_h \phi \in AP\cap D(\h)$ for each $h > 0$. By Definition 2.6, we have  $\tilde{f}'= (\phi_h-\phi)+\n \not =\n$. So the derivative depends on the choice of $f\in \tilde{f}$.

Definition 2.6  would  be correct   if only $\n \st BUC$ and $D(\h)=\{f\in BUC(\r,X): f'\in BUC(\r,X)\}$ and $Y= BUC/\n$ with $\n$ satisfying $f\i \n$ and $f'\in BUC$ implies $f'\in \n$. In this case if $\tilde{f}=\n$, then $\tilde{f}'=\n$ and so $\widetilde{\h}$ is well defined on $D(\widetilde{\h})$.

Q5. The proof of Cor 2.20 in [M] uses Theor.s  2.12 and 2.18;
         which of these (or both) is not correct?

The main mistake is the following argument in the proof of Lemma 2.7: (I rewrite from [3,  Lemma 2.7, lines -15 to -11 ])

"we have

\qquad \qquad \qquad \qquad $\la g(t)-g'(t)= h(t)$, for all $t\in \r$,

\noindent where $h$ is a function in $\n$. However, since $h\in \n$ from the condition (iii) of Condition $F$, $g$ must be in $\n$, so $\tilde {g}=0$."

\

Counter examples:  If $h= R(\la,\h) \phi$ with $\phi$ of example 1 below or $H= R(\la,\h) f$ with $f$ of example 2 below, then $h, H\in AP(\r,\cc)$ (see below) but $h'-\la h = \phi \not \in AP(\r,\cc)$  and $H'-\la H = f \not \in AP(\r,\cc)$.

 Indeed, if Re ${\la} >0 $, then by [3, (2.1)] (for the  second identity)

\qquad $h (s)= (R(\la,\h) \phi)(s) = \int_0^{\infty} e^{-\la \eta}\phi (s+\eta)\, d\eta= \check {f}_{\la}*\phi (s) $.

\noindent Here $f_{\la}(t) = \cases{ e^{-\la t},
\text{\,\, if\,\,} t \ge 0}\\ { 0,\text {\quad \,\,\, if \,} t<
 0}\endcases$,  \qquad
 $\check {f}_{\la} (t)={f}_{\la} (-t) $, $t\in \r$.

 \noindent  We have $f_{\la}\in L^1 (\r,\cc)$, $\phi\in BC(\r,\cc)$ and so $h\in BUC(\r,\cc)$. Similarly, the case Re ${\la}  <0 $. Since $sp_{AP}\phi =\emptyset$, we conclude $sp_{AP} h =\emptyset$. It follows $h\in AP$ by [1, Corollary 2.3 (A)]. In  the same way, $H\in AP (\r,\cc)$.

 \

This is why the author wrote his added in proofs.  Hypothetically, if $h', H' \in AP$, both $\phi, f$ would be in $AP$.

I think we should concentrate on refuting Lemma 2.7, since as you realized even the definition 2.8 is not correct since the function  $R(\la,\widetilde {\h})\tilde {f} =\n$ but $f\not \in \n$.

\

       Corollary  2.20 becomes true if one adds the assumption "$f$ uniformly continuous on $\r$",
       for a proof  see e.g.  [1,, p. 124, Corollary 2.3A].

4.  If one adds the condition (*) $\h(D(\h) \cap \n)  \st  \n$  (Definition 2.1 of [3]) to  Definition 2.3 of [3] as the author  suggests in his "Added to the proofs",
      then not even $\n = AP(\r,\cc)$ or $AA$ satisfy  (*):

\proclaim{Example 3} $\n\in \{AP (\r,\cc), AA (\r,\cc), REC_b(\r,\cc)\}$ does not satisfy \,
 $\h (D (\h)\cap \n)\st \n$.

(a) If $f$  is defined by $ f(t) = M_1 \phi (t) =\int_{0}^{1} \phi (t+s)\, ds$ with $\phi$ of Example 2,  then $f\in AP$ but $f'=(\phi_1-\phi)\not\in AP $.

(b) If $F= Pf$ with  $f$ of example 1, then $F\in \n$ but $f=F'\not \in \n$.
\endproclaim

    See also  the answer to Question 5 above.

In the following it is easily verified that $\n$ satisfies Condition $F$ in [3, Definition 2.3].

\proclaim{Example 4} Let $\n\in \{BUC(\r, X), BUC(\r, X)+\frak{g}\cdot BUC(\r, X)\} $, where   $X$  a Banach space and $\frak{g} (t)=e^{it^2}$. If  $f\in L^{\infty} (\r,X)$ then $sp_{\n} (f)=\emptyset$.

Indeed, $f*\phi \in BUC(\r, X)$ for each $\phi\in L^{1} (\r,\cc)$.
\endproclaim

\Refs

Here follow some questions for you (not for the editors)  :

Q1 . I do not understand Def 2.8 of [M] of $sp_{\n}$ : $\n$ does not appear in
        this definition ??
        (And with the correct definition is then always  $sp_{\n} \st sp_0
        $ ?)

Answer: $\n$ appears implicitly because the class $\tilde{f}$ is an element of $BC/\n$. But  Definition [3, 2.8] contains also a misprint, there $R(\la,\tilde \h)f$ should be replaced by $R(\la,\widetilde {\h}) \tilde{f}$. If $\n= \{0\}$ Definition 2.8 gives the Carleman spectrum of $f$ which is equal to Beurling spectrum  of $f$. In Analysis Paper 122 this is discussed for reduced spectrum. See also [3, Remark 2.9]

Q2. I do not see that the $f = \sum h_n$ of Ex.1 of [C]  is bounded ?

Example 1 is in fact  Example 2.11 of our Analysis paper 118.  $f$ is bounded because $supp\,\, h_k \cap supp\,\, h_n =\emptyset$ if $k\not = n$.

This example is very important because it shows that if $F \in AP$ and the derivative of F is continuous and bounded, then $F'$ is not even recurrent or Poisson stable. This is a missing example in the literature. It should be indicated to the  theorems: if the derivative is uniformly continuous in some ($\la$-class of almost periodicity),
then  the derivative belongs to this class.

Q3. Who showed first that  $sp_A$ f empty is equivalent with $f  \in \A$
         a) for $\A = AP$,    b)  for general $\A$ ?    (Assuming $f  \in  C_{ub}$)

For $AP$ Who showed first, I do not know But may be one can say Loomis but for general $A$, it is in my paper in [Diss. Math. 338 (1995), theorem 4.2.1)].

Q4. What is Ex. 4 of [C] for in connection with [M] ?

This example shows that $\Lambda = \emptyset$ but $\Lambda_{\n} (\Bbb{X})=\Bbb{Y}$ (see [3, (2.10) and  definition of $\Bbb{Y}$ above Definition 2.6]). This explains that the argument of the proof of Corollary 2.20 is not correct.
\enddocument